\documentclass[12pt]{article}
\usepackage{amsmath}
\usepackage{amssymb}
\usepackage[textures]{graphics}
\usepackage{color}
\usepackage{multicol}

\usepackage{graphicx}

\newtheorem{Def}{Definition}[section]
\newtheorem{Lem}[Def]{Lemma}
\newtheorem{theorem}[Def]{Theorem}
\newtheorem{Proposition}[Def]{Proposition}

\newtheorem{Rem}[Def]{Remark}
\newtheorem{example}[Def]{Example}

\def\a{{\langle\eta\rangle}}
\def\ai{{\langle\eta\rangle^{-1}}}
\def\aa{{ |\eta|}}
\def\aai{{ |\eta|^{-1}}}
\def\aaa{{ \|\eta\|}}
\def\b{{\langle\xi\rangle}}
\def\bi{{\langle\xi\rangle^{-1}}}
\def\bb{{ |\xi|}}
\def\bbi{{ |\xi|^{-1}}}
\def\bbb{{ \|\xi\|}}
\def\ab{{\xi\otimes\eta }}
\def\achl{{\hat{\triangleleft}}}

\def\acbl{{\bar{\triangleleft}}}

\def\actl{{\tilde{\triangleleft}}}
\def\actr{{\tilde{\triangleright}}}
\def\acl{{\triangleleft}}
\def\acr{{\triangleright}}

\begin{document}

{\renewcommand{\baselinestretch}{1}
   \title{\bf FURTHER RESULTS ON COSET REPRESENTATIVE CATEGORIES \\ }
     \author{ M.M\ Al-Shomrani* \& E.J.\ Beggs\dag
\vspace{0.15in}\\}
\maketitle   { }
{\center{*Department of Mathematics, Science College \\
King Abdulaziz University\\ Saudi Arabia \\ 
\dag Department of Mathematics \\ University of Wales Swansea \\ SA2 8PP U.K. \\}}
\begin{abstract}
This paper is devoted to further results on the nontrivially associated
categories
$\mathcal{C}$ and $\mathcal{D}$, which are constructed from
a choice of coset representatives for a subgroup of a finite group.
  We look at the construction of
integrals in the algebras $A$ and $D$ in the categories.  These
integrals are used to construct abstract projection operators to
show that general objects in $\mathcal{D}$ can be split into a sum
of simple objects.  The braided Hopf
algebra $D$ is
shown to be braided cocommutative, but not braided commutative.
Extensions of the categories and their connections with
conjugations and inner products are discussed.
\end{abstract}


\section{Introduction}
This paper is devoted to further results on the nontrivially associated
categories
$\mathcal{C}$ and $\mathcal{D}$, which are constructed from
a choice of coset representatives $M$ for a subgroup $G$ of a finite group
$X$
in \cite{BNT}. There are objects $A$ and $D$, in the categories
$\mathcal{C}$ and $\mathcal{D}$ respectively, which are algebras associative
in the
categories (but not associative in the `usual' sense).
The paper \cite{shombeggs} shows that the braided category
$\mathcal{D}$, which can be thought of as the double of $\mathcal{C}$, is a
modular category.
These constructions can be thought of as a nontrivially associative
version of bicrossproducts \cite{Maj1,Tak1}. For other information on braided Hopf 
algebras, see  \cite{MajBook,Tak3}.
To save a large amount of paper, we will assume the notation and results of
\cite{BNT}.

In this paper we discuss integrals for the algebras $A$ and $D$, and
projections on representations of the algebras. Then we show that $D$ is
braided cocommutative, but not braided commutative.
We discuss an additional class of morphisms (type $B$) which can be
defined
on the category $\mathcal{C}$, and their relations to the previous morphisms
(called type $A$)
and to the tensor product via a functor, Bar. Finally we discuss inner products, and their
relation to type $B$ morphisms, and give an example
involving an integral for the algebra $A$. 

Throughout the paper we assume that all groups mentioned, unless
otherwise stated,  are finite, and that all vector spaces are
finite dimensional over a field $k$, which will be denoted by
$\underline{1}$ as an object in the category.

\section{Integrals in $\mathcal{D}$}
In the literature there are two definitions of integral on a Hopf algebra,
depending
on whether it is viewed as an operator or an element:
\begin {Def}\,\,\cite{MajBook}\,\,Let $H$ be a  Hopf algebra
over the field $k$.  A left integral on $H$ is a non identically
zero linear map  \, $\int : H \rightarrow k$\, satisfying $
(\rm{id} \otimes \int )\,\circ \,\Delta= \eta \, \circ \,\int $\,.
Correspondingly right integrals are defined by 
$
(\int\otimes\rm{id}  )\,\circ \,\Delta= \eta \, \circ \,\int $\,.
If $\int 1 =1$,
then the integrals are called normalised.
\end {Def}
\begin {Def}\label{intinh}\,\cite{MajBook, Larson}\, Let $H$ be a  Hopf
algebra
over the field $k$.  A non-zero element $\Lambda \in H$ is called
a left integral if $h \Lambda  = \epsilon (h) \Lambda$ for all $h
\in H$. Similarly,  $\Lambda \in H$ is called a right integral if
$ \Lambda h = \epsilon (h) \Lambda$ for all $h \in H$. An element
$\Lambda \in H$ is called integral if it is both right and left
integral.  Integrals are normalised if $\epsilon (\Lambda)=1$\,.
\end{Def}
These definitions are of course connected. For example,
given a left integral $\int : H \rightarrow k$,
if we set $\Lambda^{*} \in H^{*}$ to be equal
to\\
\setlength{\unitlength}{0.5cm}
\begin{picture}(10.,10.)\thicklines

\put(4.,0.){\includegraphics[width=2in]{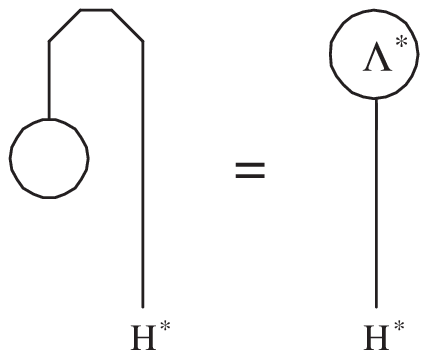}}
\put(15.5,0.){\includegraphics[width=2in]{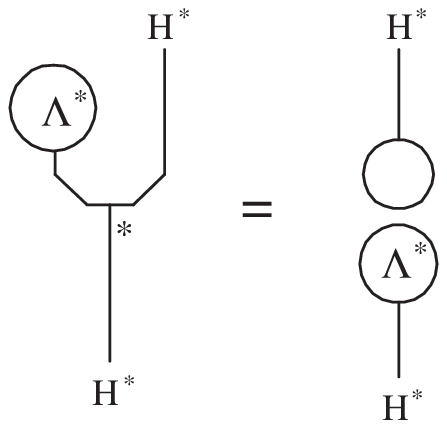}}

\put(13.5,4.3){\text{,  then}}

\put(25.5,4.8){\text{,}}

\put(5.3,4.35){\text{$\int$}}

\put(23.8,5.35){\text{$\epsilon^{*}$}}

\end{picture}
\centerline{ \rm {Figure 1}
\qquad\qquad\qquad\qquad\qquad\qquad\qquad\qquad\qquad\qquad
\qquad\qquad\qquad }
\\

\noindent i.e.\,\,$\Lambda^{*} \in H^{*}$ is a right integral in  $H^{*}$.

Now we consider our categories and give specific examples of
integrals.  First we give a definition and two useful results from 
 \cite{Maj3}. The reader will find the diagramatic proofs, which are
 quite complicated, in \cite{Maj3}. 

\newpage

$$
$$
$$
$$
\begin {Lem}\label{fig10}\ \cite{Maj3}
$$
$$
\setlength{\unitlength}{0.5cm}
\begin{picture}(10.,10.)

\put(7.9,0.){\includegraphics{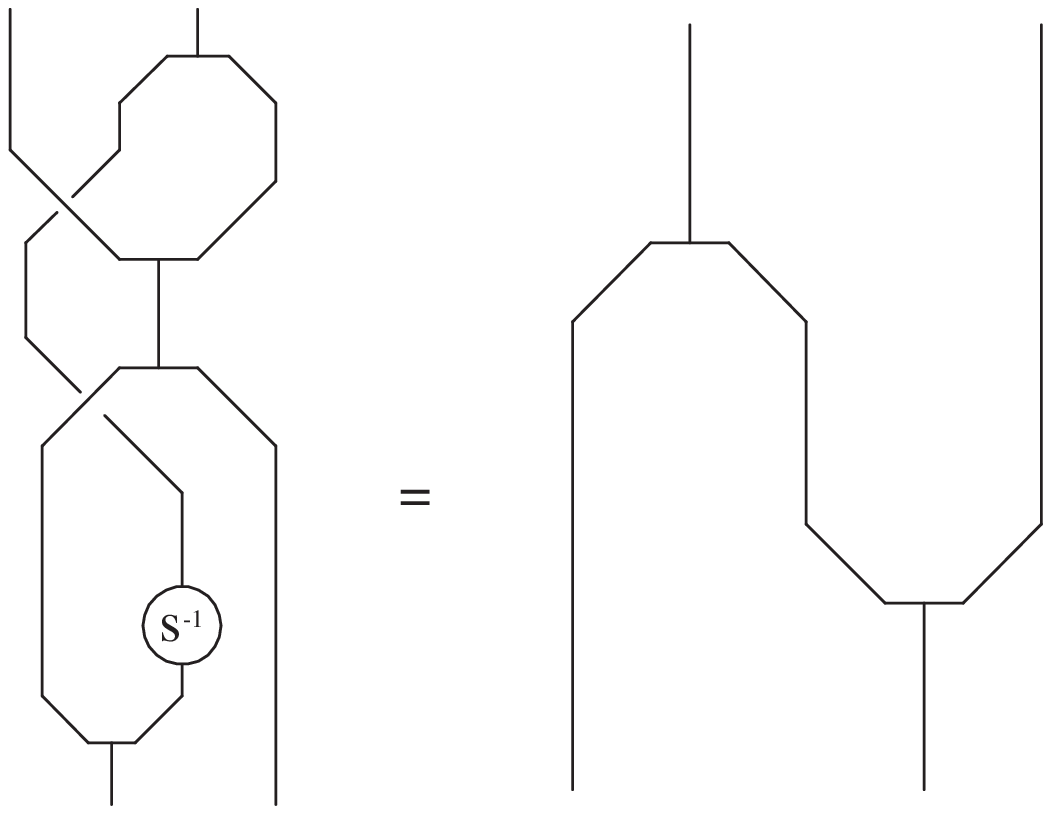}}

\end{picture}

\centerline{ \rm {\quad Figure 2 } }
\end {Lem}

\begin {Def} \cite{Maj3}\label{defoflint}\,\,\,For a braided Hopf algebra $H$, define\,
$\int : H \rightarrow k$\, by
$$
\int (h)= \rm{trace}\big(L_{h} \circ S^{2}\big) \quad \forall h
\in H \,,
$$
where $L_{h}$ is the left multiplication by $h$.  This can be illustrated by
the following diagram:\\
\vspace{1.5cm}

\begin{picture}(10.,10.)

\put(10.5,0.){\includegraphics{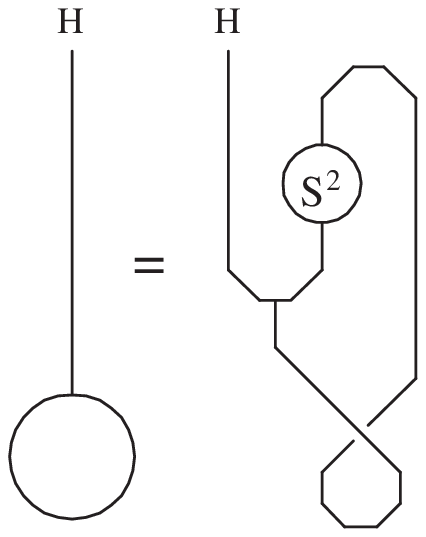}}

\put(11.8,2){$\mathbf{\int}$}

\end{picture}
\centerline{ \rm { Figure 3}\qquad\quad\qquad\qquad\qquad\qquad\qquad\qquad
\qquad\qquad\qquad\qquad\qquad\qquad}
\end {Def}

\newpage

\vspace{0.5cm}
\begin{Proposition}\label{fig15} \cite{Maj3} The map $\int$
defined in  \ref{defoflint} is a left integral, i.e.
$$
$$
\begin{picture}(10.,10.)

\put(9.,0.){\includegraphics{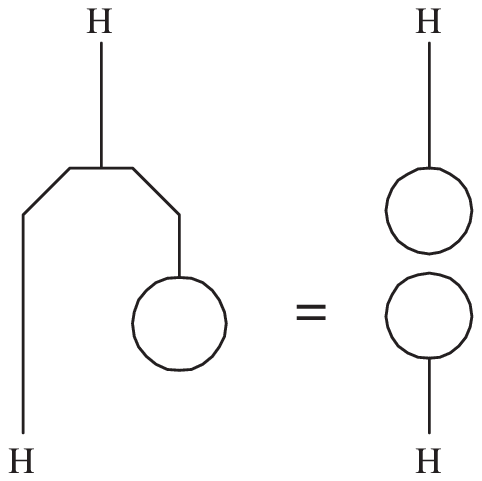}}

\put(13.5,3.3){$\int$}

\put(18.55,5.7){$\int$}

\put(18.6,3.75){$\eta$}
\end{picture}
\centerline{ \rm {Figure 4 \qquad}
\qquad\qquad\qquad\qquad\qquad\qquad\qquad\qquad\qquad\qquad\qquad}
\end{Proposition}

\begin{Proposition}\label{roint}\,\,\,In the braided tensor category
$\mathcal{D}$,
for an element $\rho \in D$ satisfying
$$
$$
$$
\setlength{\unitlength}{0.5cm}
\put(-4.0,0.0){\begin{picture}(10.,10.)
\put(-8.8,-0.6){\includegraphics[scale=1.2]{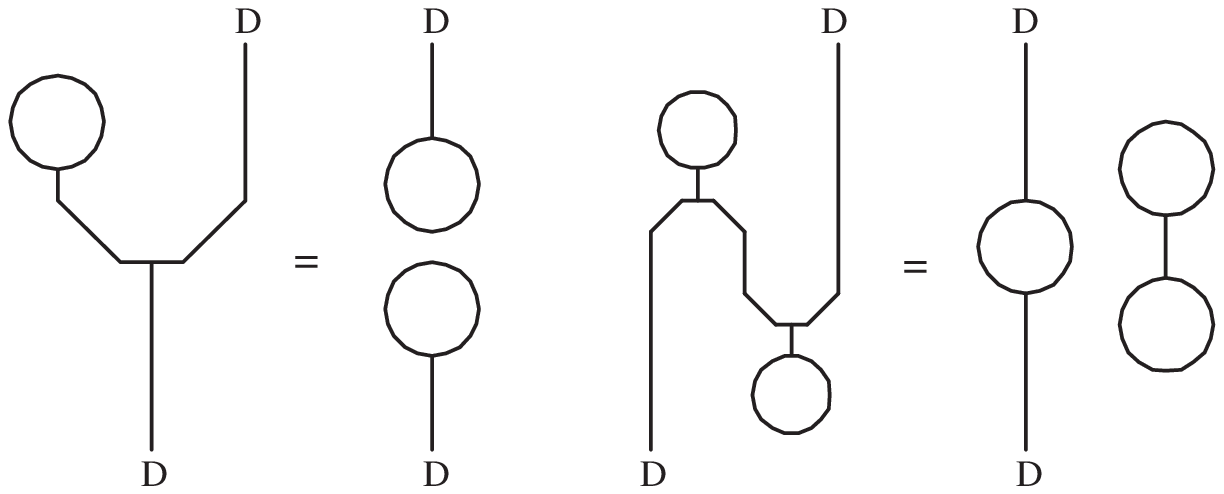}}
\put(3.3,5.2){\text{\rm{, we have}}} \put(21.5,5.2){\text{\rm{.}}}
\put(-7.2,8.6){\text{$\rho$}}\put(1.8,4.2){\text{$\rho$}}
\put(1.8,7.1){\text{$\epsilon$}}
\put(10.65,2.0){\text{$\int$}}\put(8.3,8.6){\text{$\rho$}}
\put(16.0,5.4){\text{$S^{-1}$}} \put(19.8,7.6){\text{$\rho$}}
\put(19.8,3.6){\text{$\int$}}
\end{picture}}\qquad\qquad\qquad
$$
\end{Proposition}
\textbf{Proof.}
$$
$$
$$
$$
$$
$$
$$
$$
$$
$$
$$
$$
$$
$$
$$
\setlength{\unitlength}{0.5cm}
\begin{picture}(10.,10.)
\put(-9.0,-01.5){\includegraphics[scale=1.2]{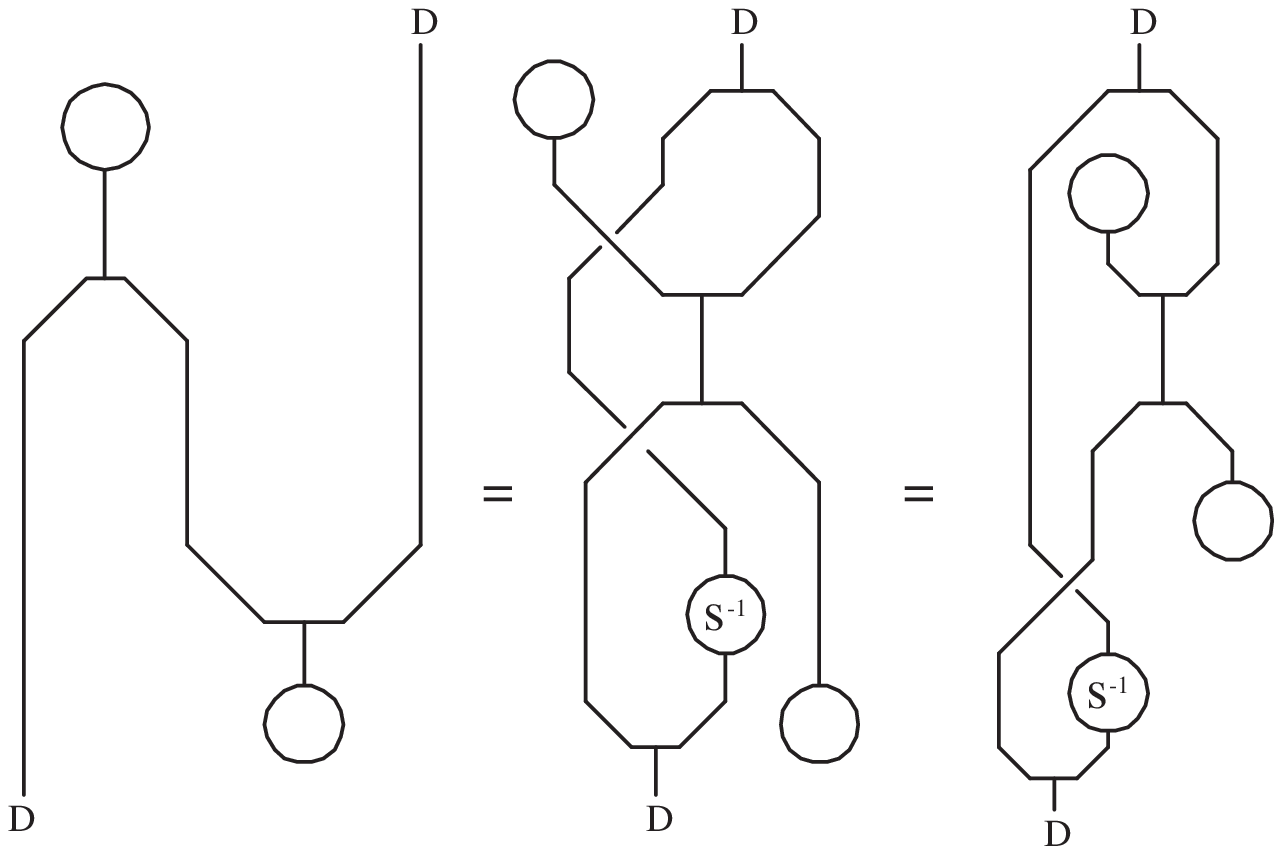}}
\put(-5.2,16.5){\text{$\rho$}}\put(5.5,17){\text{$\rho$}}
\put(19.1,14.8){\text{$\rho$}} \put(12.05,1.7){\text{$\int$}}
\put(-0.55,1.7){\text{$\int$}}\put(22.05,6.7){\text{$\int$}}

\end{picture}\qquad\qquad\qquad
$$
$$
$$
$$
$$
$$
$$
$$
$$
$$
$$
$$
$$
$$
\setlength{\unitlength}{0.5cm} \put(-4.,0){\begin{picture}(10.,10.)
\put(-9.0,-0.5){\includegraphics[scale=1.2]{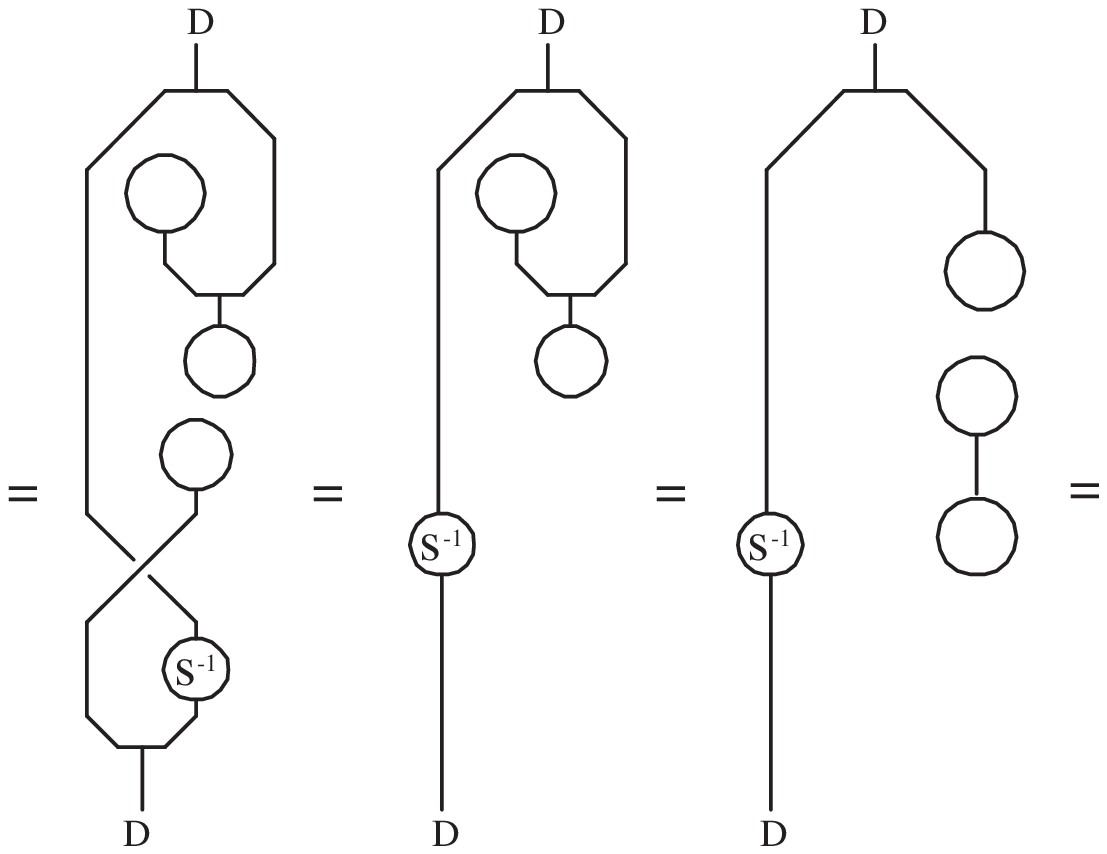}}
\put(-4.2,15.8){\text{$\rho$}}\put(4.4,15.8){\text{$\rho$}}
\put(-3,11.5){\text{$\int$}}\put(-3.5,9.3){\text{$\eta$}}
\put(5.5,11.5){\text{$\int$}} \put(15.5,7.2){\text{$\int$}}
\put(15.5,10.7){\text{$\rho$}} \put(15.7,13.8){\text{$\epsilon$}}
\put(19.5,8.5){\text{R.H.S.}}
\end{picture}}\qquad\qquad\qquad
$$
Note that we have used lemma $\ref{fig10}$ for the first equality
and proposition $\ref{fig15}$ for the third equality.  The
hypothesis for $\rho$ has been used  for the fifth
equality.$\quad\square$

\medskip

The definition of integral does not require the category to be
braided.  Here we give an example of an integral in $\mathcal{C}$.
\begin{Proposition}\,
Let $A$ be the   algebra  in the category $ \mathcal{C}$ defined
in \cite{BNT}, then the element $\rho= \sum_{u\in G}\delta_{e} \otimes
u$ for $u \in G$ is an integral element.

\end{Proposition}
\textbf{Proof.} \,\,\,We need to prove that $\rho=
{\sum}_{u}\delta_{e} \otimes u$ is both right and left integral,
so for any element $h=(\delta_{t} \otimes v )\in H$  we have
\begin{equation*}
\begin{split}
\rho h =(\sum_{u}\delta_{e} \otimes u)(\delta_{t} \otimes v )
=\sum_{u}(\delta_{e} \otimes u)(\delta_{t} \otimes v )
=\sum_{u}\delta_{t,e \acl u}\,\delta_{e \acl \tau (a, b)} \otimes
\tau (a, b)^{-1} uv,
\end{split}
\end{equation*}
where $a= \langle \delta_{e} \otimes u  \rangle$ and $b= \langle
\delta_{t} \otimes v  \rangle$.  But we know that $e \acl u = e$
and $e \acl \tau (a, b)=e$.  Moreover, $e \cdot \langle \delta_{e}
\otimes u  \rangle = e \acl u =e$. Also because $a= \langle
\delta_{e} \otimes u  \rangle=e$, then $\tau (a, b)=\tau (e,
b)=e$. Now  as $uv$ is an element in $G$, then we get
\begin{equation*}
\begin{split}
\rho h =\delta_{t,e } \sum_{u}\,\delta_{e } \otimes  uv
=\delta_{t,e } \rho = \epsilon (h) \rho,
\end{split}
\end{equation*}
so $\rho= \sum_{u}\delta_{e} \otimes u$ is a right integral . Next
we want to show that it is also a left integral, so we start with
\begin{equation*}
\begin{split}
h \rho  =(\delta_{t} \otimes v )(\sum_{u}\delta_{e} \otimes u)
=\sum_{u}(\delta_{t} \otimes v )(\delta_{e} \otimes u)
=\sum_{u}\delta_{e,t \acl v}\,\delta_{t \acl \tau (b, a)} \otimes
\tau (b, a)^{-1} vu,
\end{split}
\end{equation*}
where $b= \langle \delta_{t} \otimes v  \rangle$ and $a= \langle
\delta_{e} \otimes u  \rangle$.  But we know that  $t \acl e = t$.
Moreover, $\delta_{e,t \acl v}=1$ implies $e=t \acl v$ or $e \acl
v^{-1}=e=t$. Also because $a= \langle \delta_{e} \otimes u
\rangle=e$, then $\tau (b, a)=\tau (b,e)=e$. Now  as $vu$ is an
element in $G$, then we get
\begin{equation*}
\begin{split}
h \rho  =\delta_{e,t \acl v}\, \sum_{u}\,\delta_{t } \otimes vu
=\delta_{e,t } \sum_{u}\,\delta_{e } \otimes  vu =\delta_{t,e
}\,\rho =\epsilon (h) \rho. \qquad\square
\end{split}
\end{equation*}

\section{\!\!\!Projections on representations in $\mathcal{D}$
using integrals} Before going further, we recall some concepts and
results from  finite group representations. Later we will apply these to 
the braided Hopf algebra $D$ in the category $\mathcal{D}$. 

Let $V$ be a vector space, and let $W$ and $W_{o}$ be two
subspaces of $V$. Then for the direct sum $V= W \oplus W_{o}$,
$W_{o}$ is called a complement of $W$ in $V$.   The map $p$
which sends each $x \in V$ to its component $w \in W$ is called
the projection of $V$ onto $W$ associated with the
decomposition $V= W \oplus W_{o}$\,. The image of $p$ is $W$, and
$p(x)=x$ for all $x \in W$.

\begin{theorem} \,\,\cite{Serre}\, Let $\rho$ be a linear representation of
a finite group $G$
in $V$ and let $W$ be a vector subspace of $V$ stable under $G$.
Then there exists a complement $W_{o}$ of $W$ in $V$ which is
stable under G.

\end{theorem}
\textbf{Proof.}\hspace{0.5cm} Let $W_{o}$ be an arbitrary
complement of $W$ in $V$, and let $p$ be the corresponding
projection of $V$ onto $W$. We know that from the definition of
the average $p_{_o}$ of the conjugates of $p$ by the elements of
$G$:
$$
p_{_o}= \frac{1}{n}\,\sum_{t \in G} \rho_{t} \cdot p \cdot
{\rho_{t}}^{-1},
$$
where $n$ is the order of $G$.  Since $p$ maps $V$ into $W$ and
$\rho_{t}$\, preserves $W$ we see that $p_{_o}$ maps $V$ into $W$.
We have ${\rho_{t}}^{-1} (x)\in W$, If $x \in W$ hence
$$
p \cdot {\rho_{t}}^{-1} (x)= {\rho_{t}}^{-1} (x), \qquad \rho_{t}
\cdot p \cdot {\rho_{t}}^{-1}(x)=x, \qquad \text{and}\qquad
p_{_o}(x)=x.
$$
Thus $p_{_o}$ is a projection of $V$ onto $W$, corresponding to
some complement $W_{o}$ of $W$.  Moreover, we have $\,\,
\rho_{t}\, p_{_o}=p_{_o}\,\rho_{t} \quad \text{for all}\,\, s \in
G. $ If we compute $\rho_{s} \cdot p_{_o} \cdot {\rho_{s}}^{-1}$,
we find:
$$
\rho_{s} \cdot p_{_o} \cdot {\rho_{s}}^{-1}=\frac{1}{n}\,\sum_{t
\in G} \rho_{s} \cdot \rho_{t} \cdot p \cdot {\rho_{t}}^{-1} \cdot
{\rho_{s}}^{-1} =\frac{1}{n}\,\sum_{t \in G}\rho_{st} \cdot p
\cdot {\rho_{st}}^{-1}=p_{_o}.
$$
Now for $x \in W_{o}$ and $s \in G$, we have $p_{_o}(x)=0$ which
implies that
$$
p_{_o} \cdot {\rho_{s}} (x)={\rho_{s}} \cdot p_{_o} (x)=0,
$$
that is $({\rho_{s}} (x)) \in W_{o}$, which shows that $W_{o}$ is
stable under $G$. \qquad $\square$

\medskip
We return now to the right representation of the Hopf algebra $D$
in the braided category $\mathcal{D}$ supposing that $\Lambda \in
D$
is a right integral, i.e.\\
\vspace{-1.3cm}
\begin{picture}(10.,10.)

\put(10.,0.){\includegraphics{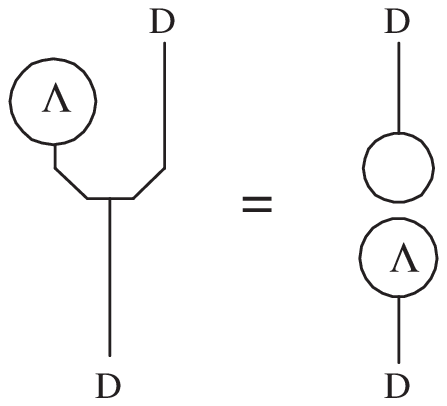}}

\put(18.5,5.1){\text{$\epsilon$}}

\end{picture}\\
\centerline{ \rm {Figure 5} }
$$
$$
\begin {Lem}\label{fig110}\

\begin{picture}(10.,10.)
\put(5.,0.){\includegraphics{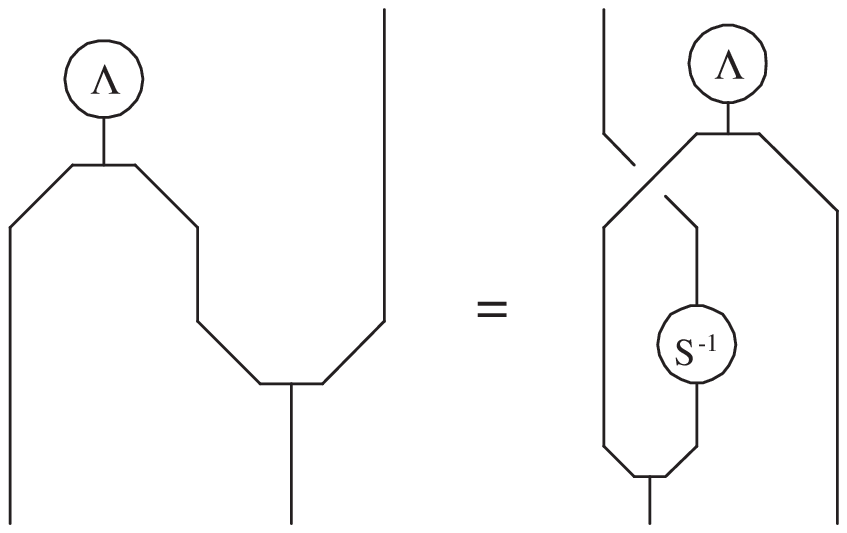}}
\end{picture}
\centerline{ \rm {Figure
6\qquad}\qquad\qquad\qquad\qquad\qquad\qquad\qquad\qquad\qquad\qquad
\qquad\qquad\quad }
\end {Lem}

\noindent
\textbf{Proof.} Using Lemma \ref{fig10}:\\
$$
$$
$$
$$
$$
$$
$$
$$
\put(0.,-9.){\begin{picture}(10.,10.)
\put(0.,0.){\includegraphics{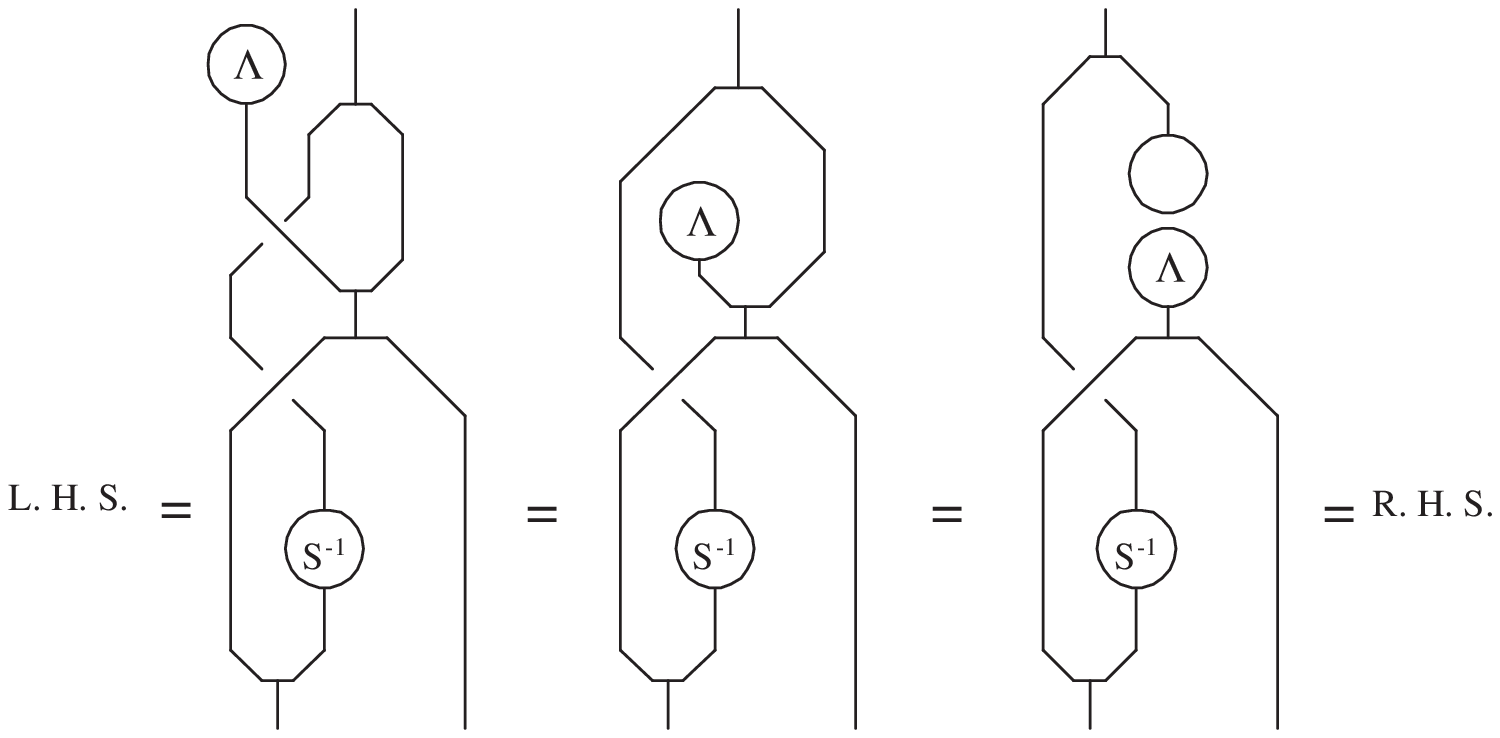}}
\put(24.15,11.4){\text{$\epsilon$}}
\end{picture}}
\begin {Def}\label{defofprot}\,\,\,For the right representations $V$ and $U$
of $D$, and a
linear map (not necessarly a morphism) $t : V \rightarrow
U$, we define $t_{_o} : V \rightarrow U $ by\\
$$
$$
$$
$$
$$
$$
\begin{picture}(10.,10.)
\put(10.,0.){\includegraphics{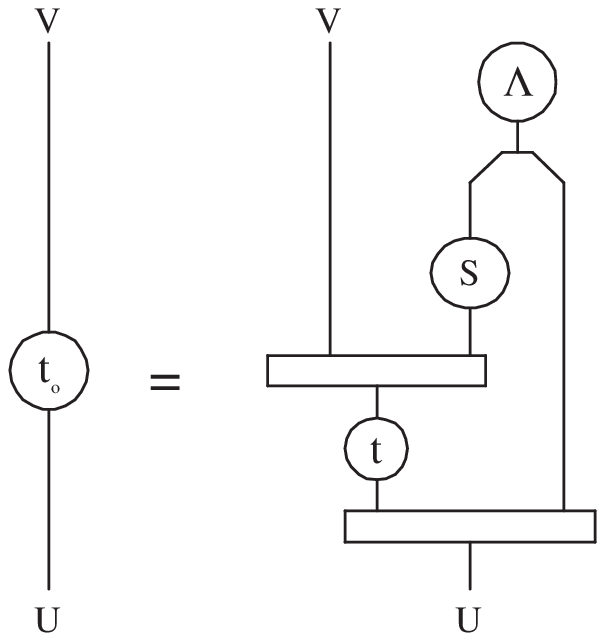}}
\end{picture}
\centerline{ \rm {\qquad Figure 7}
\qquad\qquad\qquad\qquad\quad\qquad\qquad\qquad\qquad\qquad\qquad
\qquad\qquad\qquad}
\end {Def}

\begin{Proposition}\,\label{fig12}\
The  map $ t_{_o} $ is a morphism in the category $\mathcal{D}$,
i.e.\\
\begin{picture}(10.,10.)
\put(10.,0.){\includegraphics{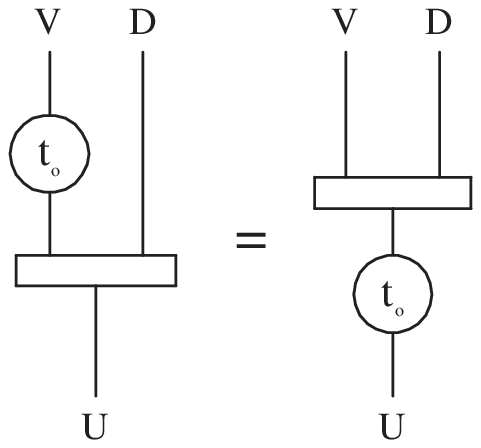}}
\centerline{ \rm {\,\,\quad Figure 8} }
\end{picture}
\end{Proposition}
\textbf{Proof.}\\
\begin{picture}(10.,10.)
\put(0.,-9.5){\includegraphics{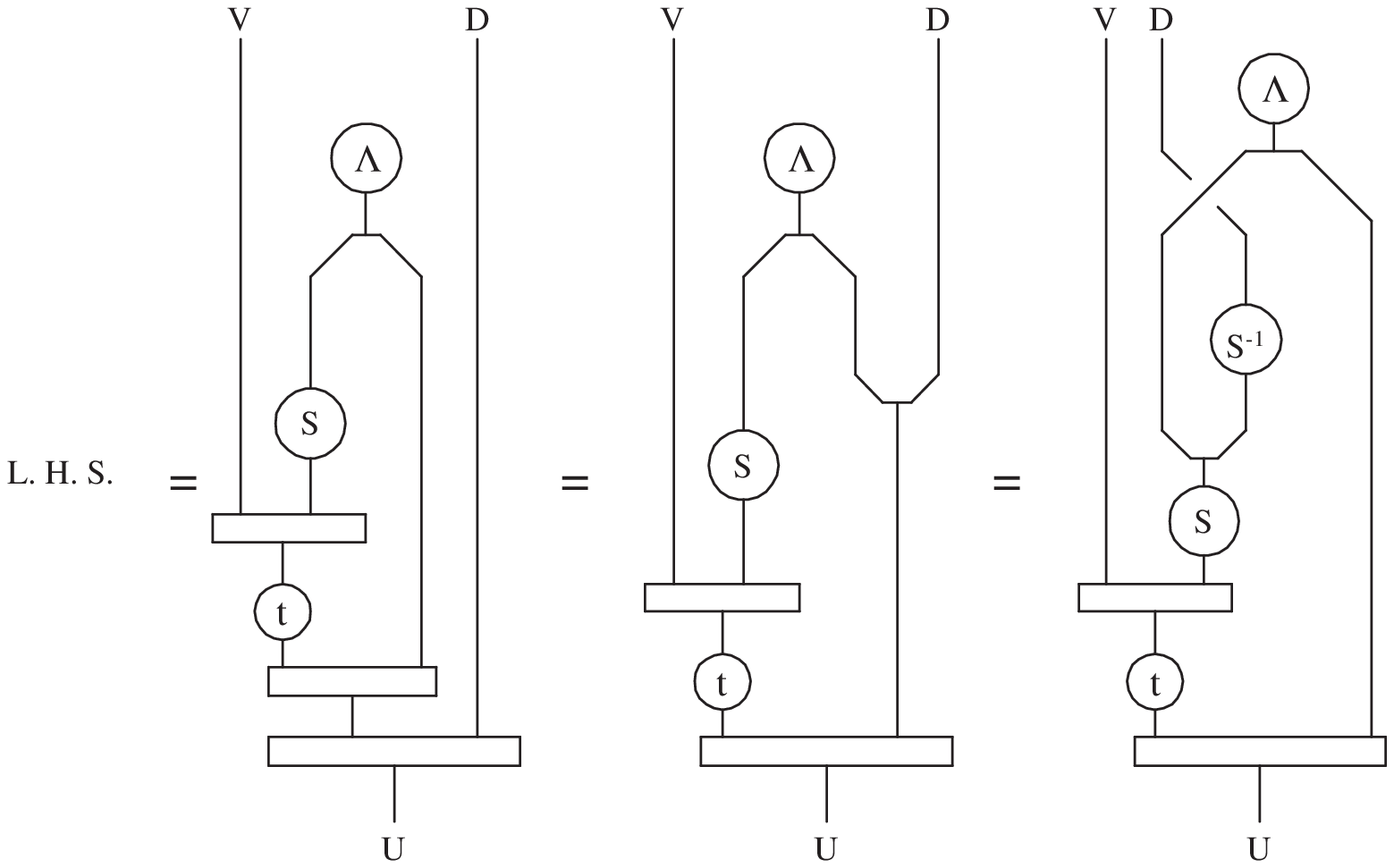}}
\end{picture}

\newpage

$$
$$
$$
$$
$$
$$
$$
$$
$$
$$
\begin{picture}(10.,10.)
\put(0.,-1.5){\includegraphics{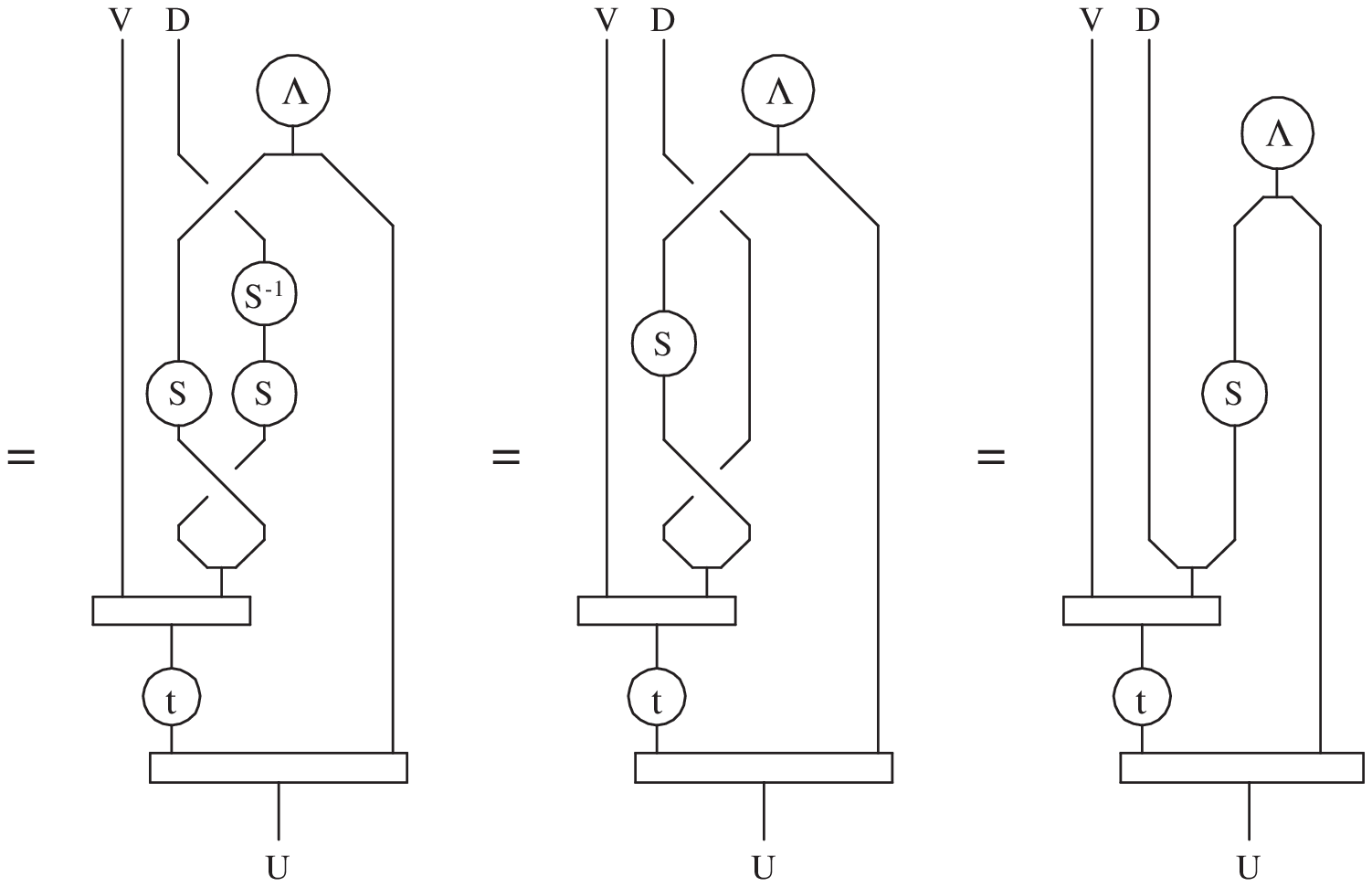}}
\end{picture}
$$
$$
$$
$$
$$
$$
$$
$$
$$
$$
$$
$$
$$
$$
\begin{picture}(10.,10.)
\put(0.,0.){\includegraphics{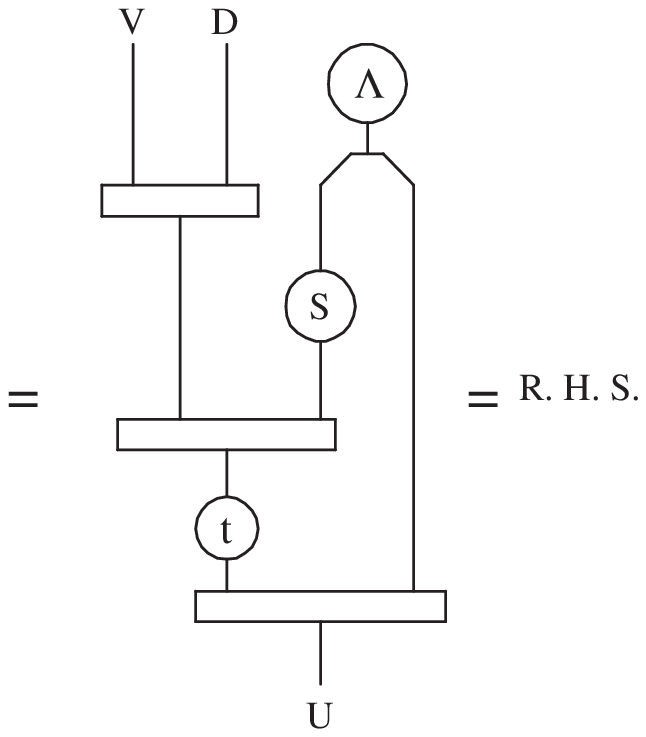}}
\end{picture}

\newpage

\begin{Proposition}\,Suppose $\epsilon (\Lambda)=1$.  Let $V$ be a
right representation of $D$, and $W \subset V$ be a
subrepresentation.  Then there is a complement $W_{o}$ of $W$ which
is also a right representation of $D$.
\end{Proposition}
\textbf{Proof.} Take any projection \,$p : V \rightarrow V$  with
image $W$.  By \ref{defofprot} we also get a morphism  $p_{_o} : V
\rightarrow V$.   Then the proof is given as follows:\\{\bf
a)}\,\,\,Show that $p_{_o}|_{_W}$ is the identity.
$$
$$
$$
$$
$$
$$
$$
$$
$$
$$
\begin{picture}(10.,10.)
\put(0.,0.){\includegraphics{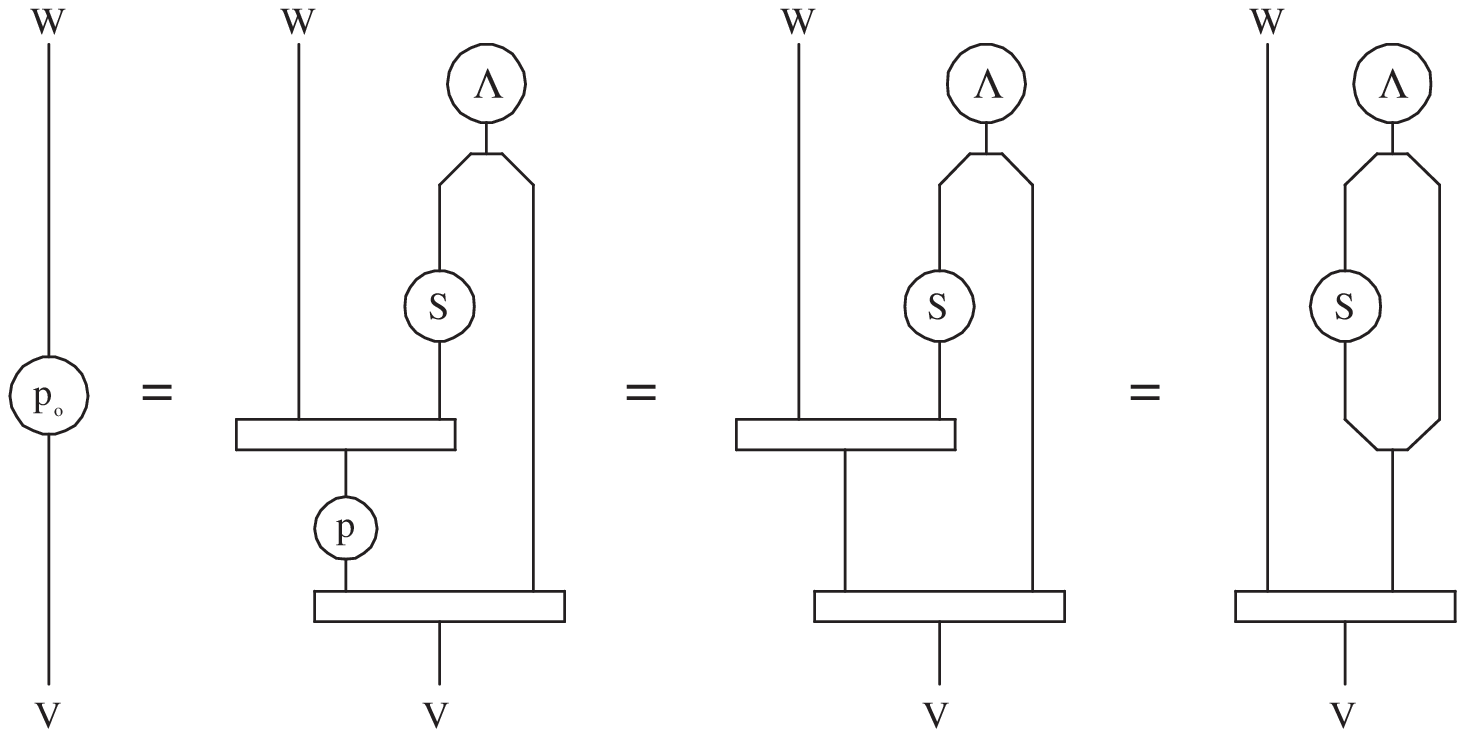}}
\end{picture}
$$
$$
$$
$$
$$
$$
\qquad\qquad\begin{picture}(10.,10.)

\put(0.,-1.){\includegraphics{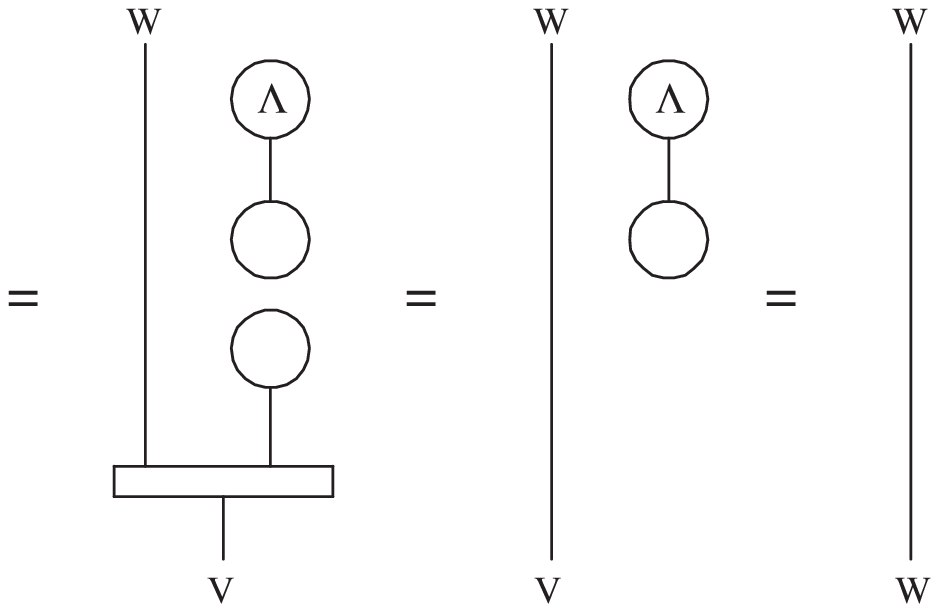}}

\put(6.15,4.7){\text{$\eta$}}

\put(6.2,7.){\text{$\epsilon$}}

\put(14.3,7.){\text{$\epsilon$}}

\end{picture}

\newpage\noindent
\textbf{b)}\,\,\,Show that the image of $p_{_o} : V \rightarrow V$ is
contained in $W$.\\
$$
$$
$$
$$
$$
$$
$$
$$
\begin{picture}(10.,10.)
\put(8.,0.){\includegraphics{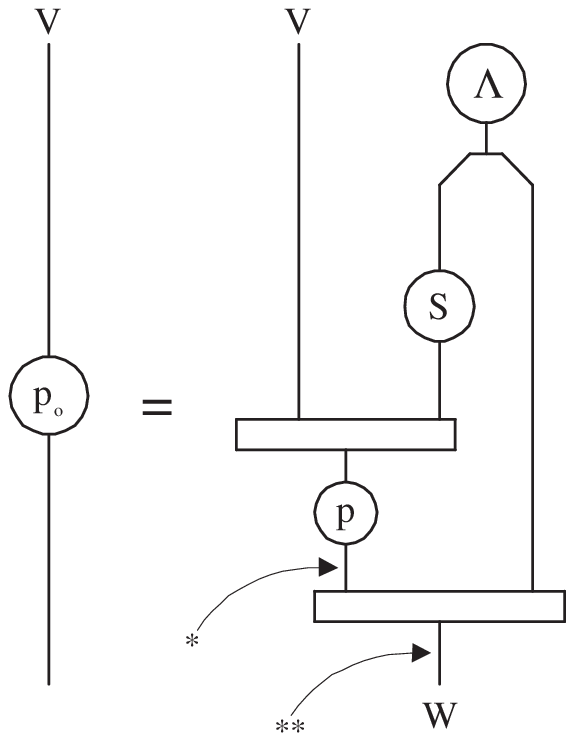}}
\end{picture}

As $p(V) \subset W$, the elements in the diagram at position $*$
is in $W$.  But as $W$ is a subrepresentation of $V$, the output
at $**$ is also in $W$.

Combining a) and b) shows that   $ p_{_o} $ is a projection.

\vspace{0.5cm} 
\noindent\textbf{c)}\,\, Show that  $W_{o}$ =ker $p_{_o}$ is
a subrepresentation.

\begin{picture}(10.,10.)
\put(0.,0.){\includegraphics{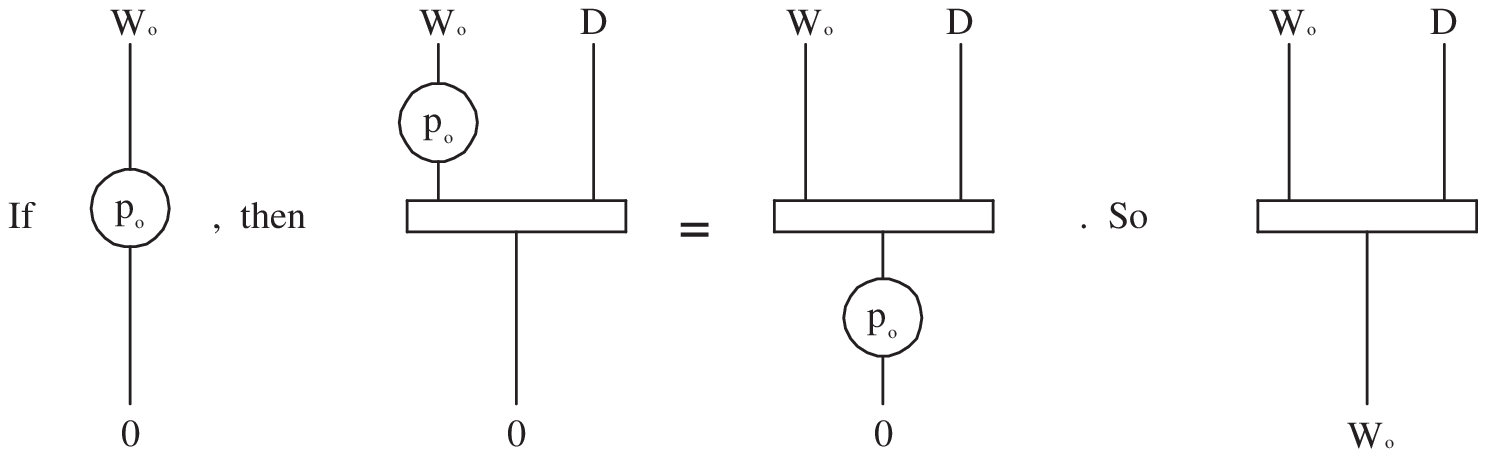}}
\end{picture}

\noindent i.e. $W_{o}$ is a right representation of $D$.

\newpage
\begin{Proposition}\,Suppose $\epsilon (\Lambda)=1$.  Let $V$ and $W$ be
two right irreducible representations of $D$. For a linear map $t
: V \rightarrow W$, by Schur's Lemma we have $t_{_o}=0$ if $V$ is
not isomorphic to $W$, and if $V=W$ then $t_{_o}= c
\,\,\rm{id}_{V}$\,.  The value of $c$ is given by \,
$c=\frac{\rm{trace(t)}}{\rm{dim}_{V}}$\,.

\end{Proposition}
\textbf{Proof.}
$$
$$
$$
$$
$$
$$
$$
$$
$$
$$
$$
$$
$$
$$
\begin{picture}(10.,10.)
\put(4.,0.){\includegraphics{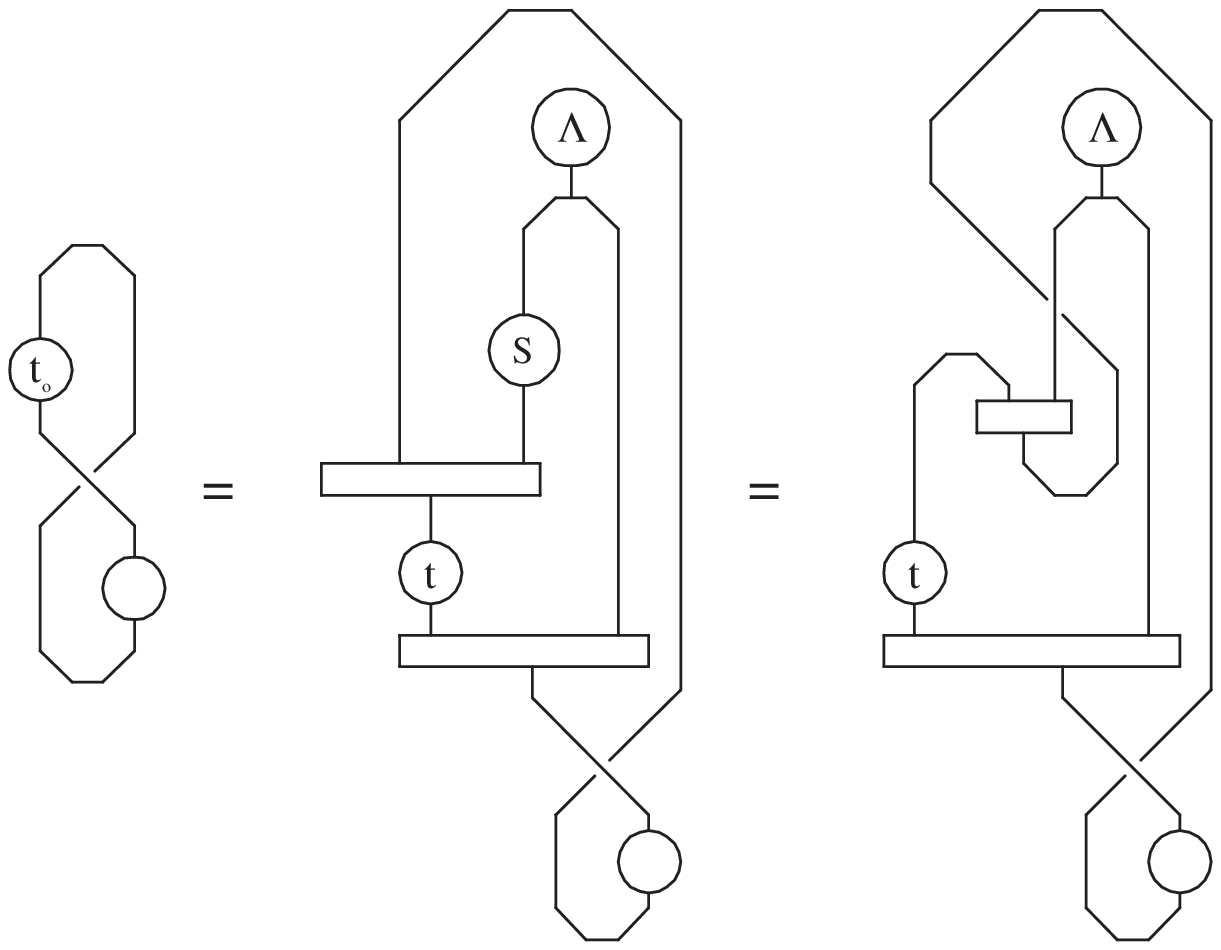}}
\put(0.15,9.7){\text{$c \, \rm{dim (V)}=$}}
\put(7.05,7.7){\scriptsize{\text{${\theta^{-1}}$}}}
\put(17.55,2.2){\scriptsize{\text{${\theta^{-1}}$}}}
\put(28.3,2.2){\scriptsize{\text{${\theta^{-1}}$}}}
\end{picture}
$$
$$
$$
$$
$$
$$
$$
$$
$$
$$
\begin{picture}(10.,10.)
\put(0.,0.){\includegraphics{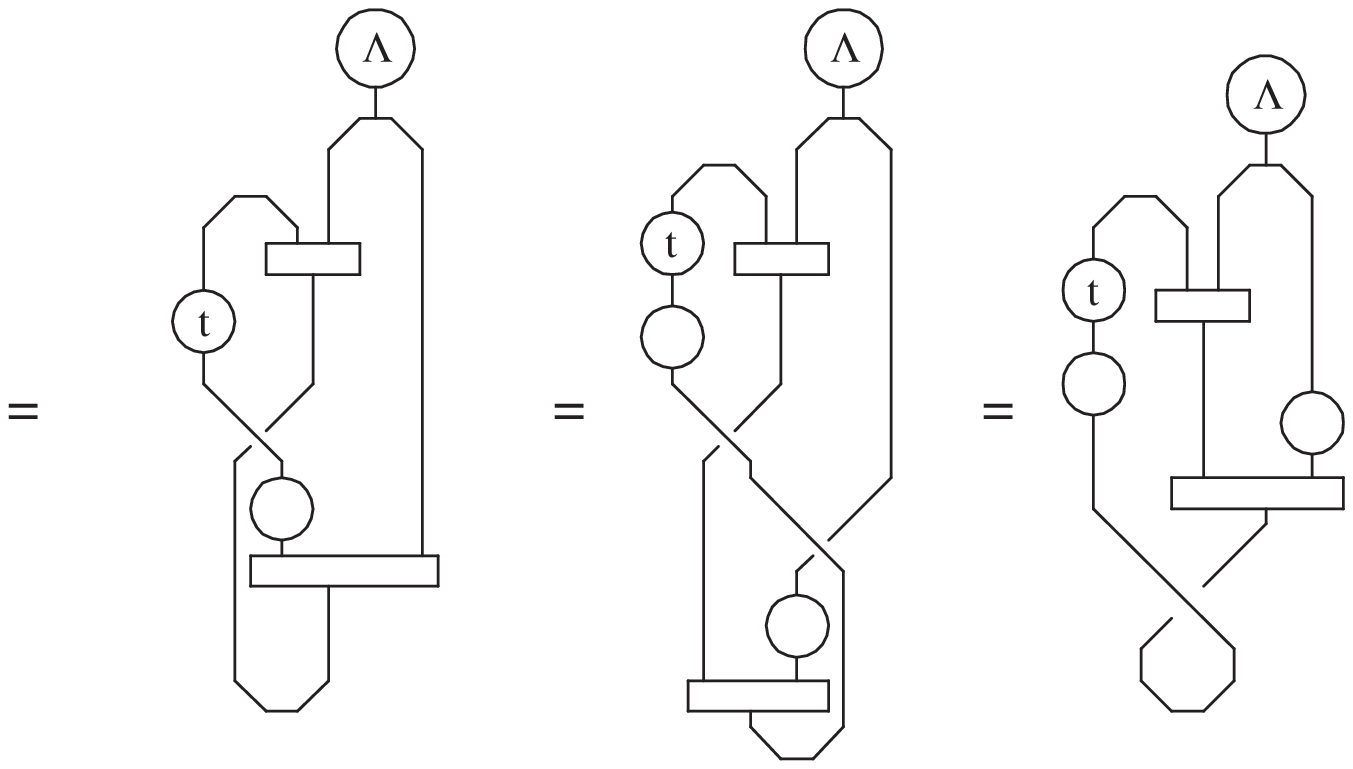}}
\put(6.45,5.63){\scriptsize{\text{${\theta^{-1}}$}}}
\put(14.35,9.13){\scriptsize{\text{${\theta^{-1}}$}}}
\put(17.05,3.23){\text{$S$}}
\put(22.9,8.18){\scriptsize{\text{${\theta^{-1}}$}}}
\put(27.5,7.33){\text{$S$}}
\end{picture}
$$
$$
$$
$$
$$
$$
$$
$$
$$
$$
\begin{picture}(10.,10.)
\put(0.,0.){\includegraphics{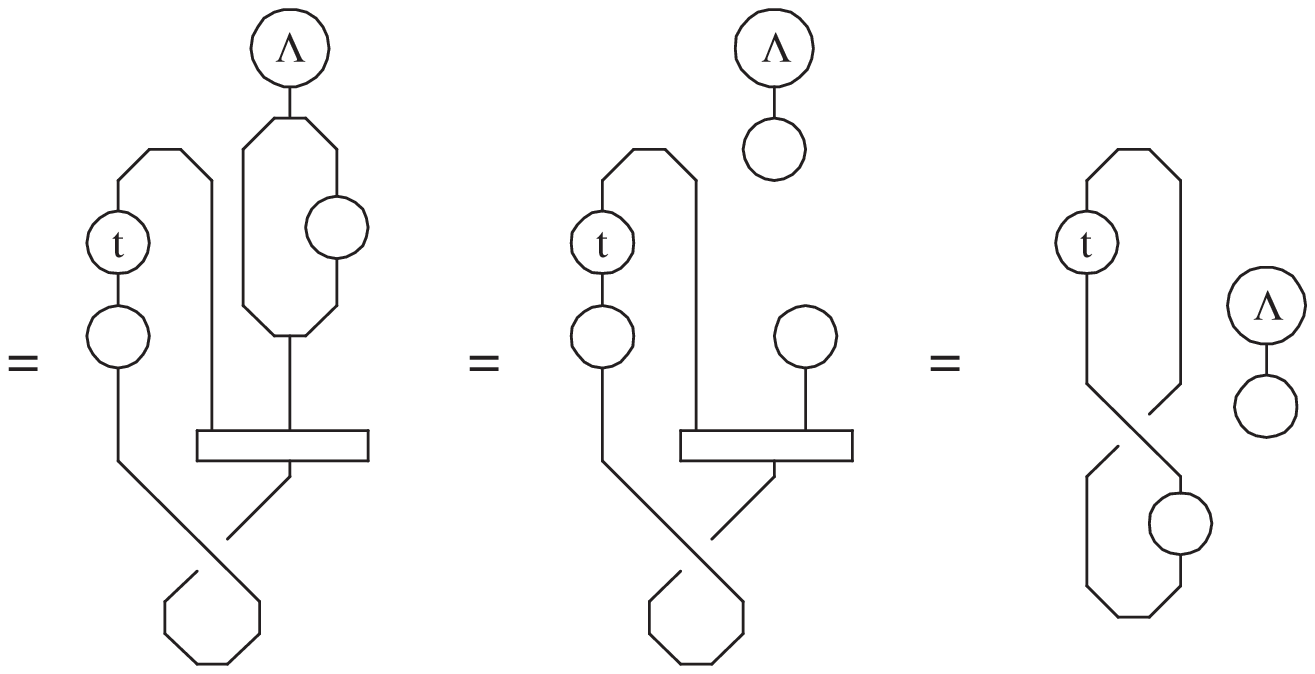}}
\put(2.8,7.18){\scriptsize{\text{${\theta^{-1}}$}}}
\put(7.3,9.38){\text{$S$}}
\put(12.6,7.18){\scriptsize{\text{${\theta^{-1}}$}}}
\put(16.3,11.1){\text{$\epsilon$}} \put(16.9,7.3){\text{$\eta$}}
\put(24.35,3.38){\scriptsize{\text{${\theta^{-1}}$}}}
\put(26.3,5.9){\text{$\epsilon$}}
\end{picture}

$= \rm{trace}(t)\,\epsilon (\Lambda)=\rm{trace}(t)\,. \quad
\square$

\section{The Hopf algebra $D$ is braided cocommutative}

We consider a braided Hopf algebra $E$ in a braided category
\,$\mathcal{S}$, in which $E$ has a right action on the objects in
\,$\mathcal{S}$ given by the morphism
$$
$$
$$
$$
\begin{picture}(6.,5.)
\put(12.5,1.){\includegraphics{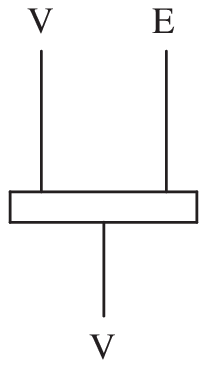}}
\centerline{\rm{\quad Figure 9} }
\end{picture}

and the action on tensor product is given by

\begin{picture}(10.,10.)
\put(7.5,0.){\includegraphics{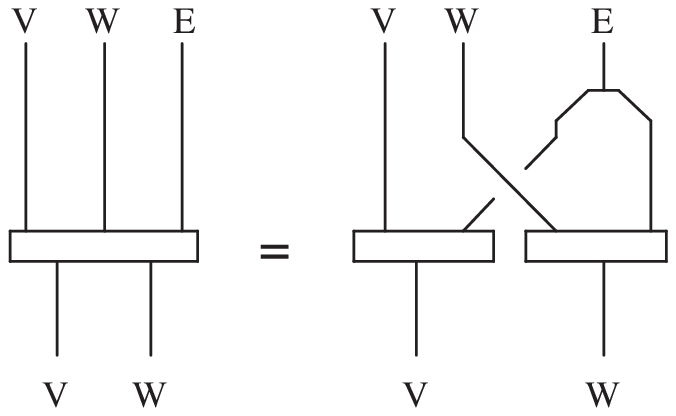}}
\centerline{\rm{ Figure 10}\quad \qquad}
\end{picture}

\begin {Def}\label{opcopr}\,\,\, The opposite coproduct, $\Delta^{op} $, for
the
algebra $E$ in  ${\cal S}$ can be  defined by the following
diagram for the representations $V$ and $W$ of $E$ in ${\cal S}$:\\
$$
$$
\vspace{0.5cm}
%
\begin{picture}(10.,10.)
%

\put(6.5,0.){\includegraphics{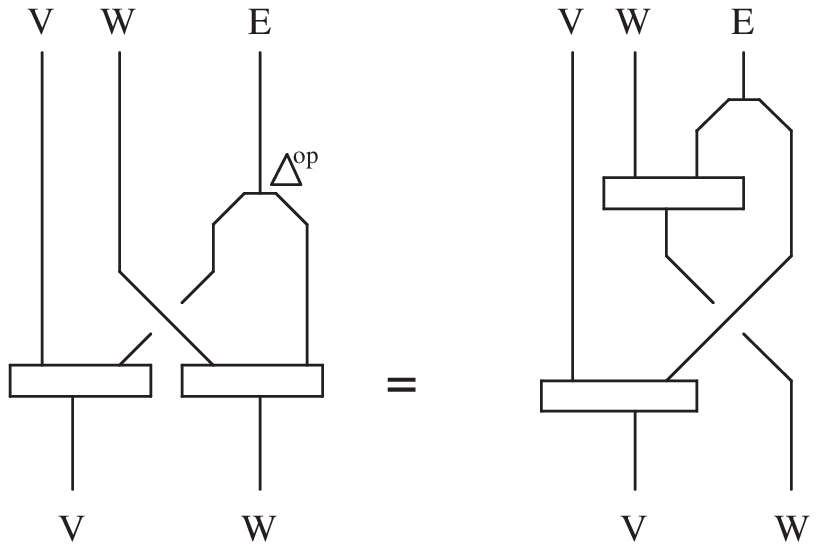}}
\centerline{\rm{ Figure 11}}
\end{picture}

\end {Def}


\vspace{-1.0cm}
\begin {Lem}\label{coopl}\
For the representations $V$ and $W$ of $E$ in ${\cal S}$, the
opposite coproduct, $\Delta^{op} $, satisfies the following
$$
$$
$$
$$
\begin{picture}(10.,10.)
\put(7.5,0.){\includegraphics{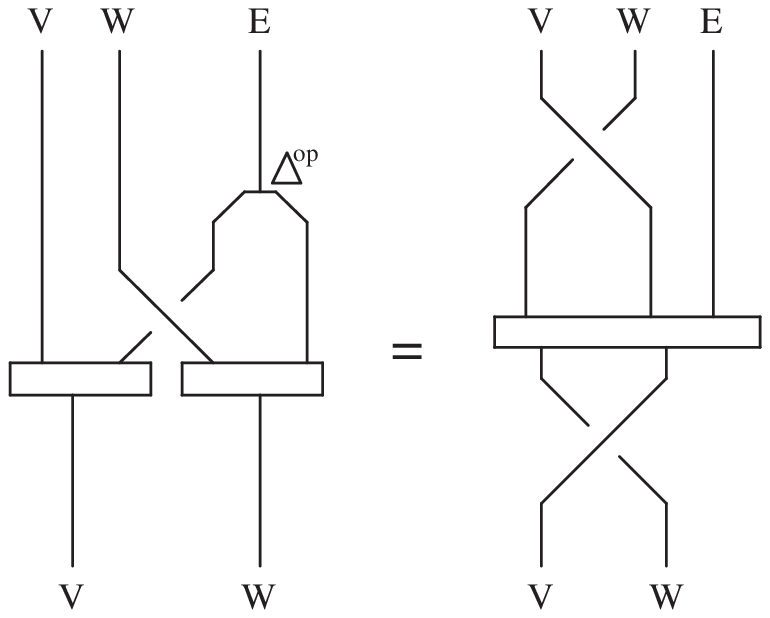}}
\centerline{\rm{\qquad Figure 12}}
\end{picture}\
\end {Lem}
\text{\bf{Proof.}}
$$
$$
$$
$$
$$
$$
\begin{picture}(10.,10.)
\put(1.,0.){\includegraphics{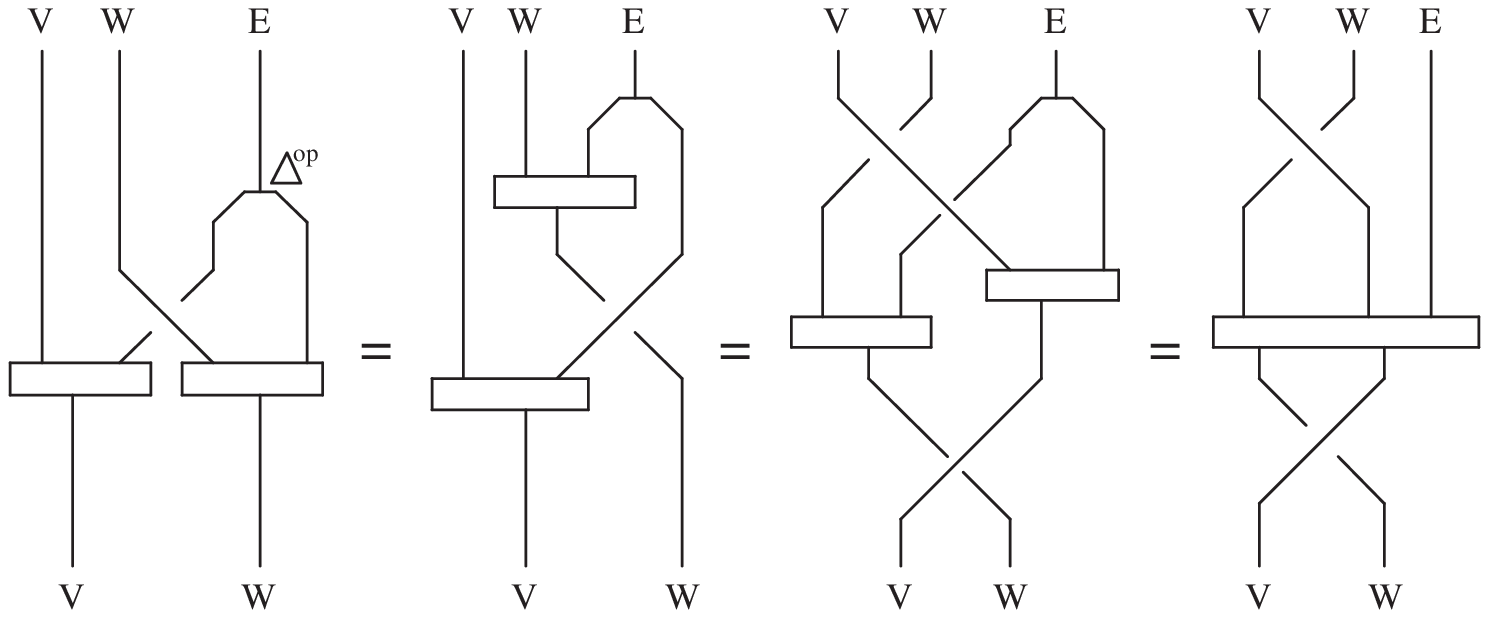}}
\end{picture}

\newpage
\begin{Proposition}\,For the algebra $E$ the opposite coproduct,
$\Delta^{op}
$,is coassociative, i.e.

$$
$$

%
\begin{picture}(10.,10.)
%
\put(2.,-1.){\includegraphics{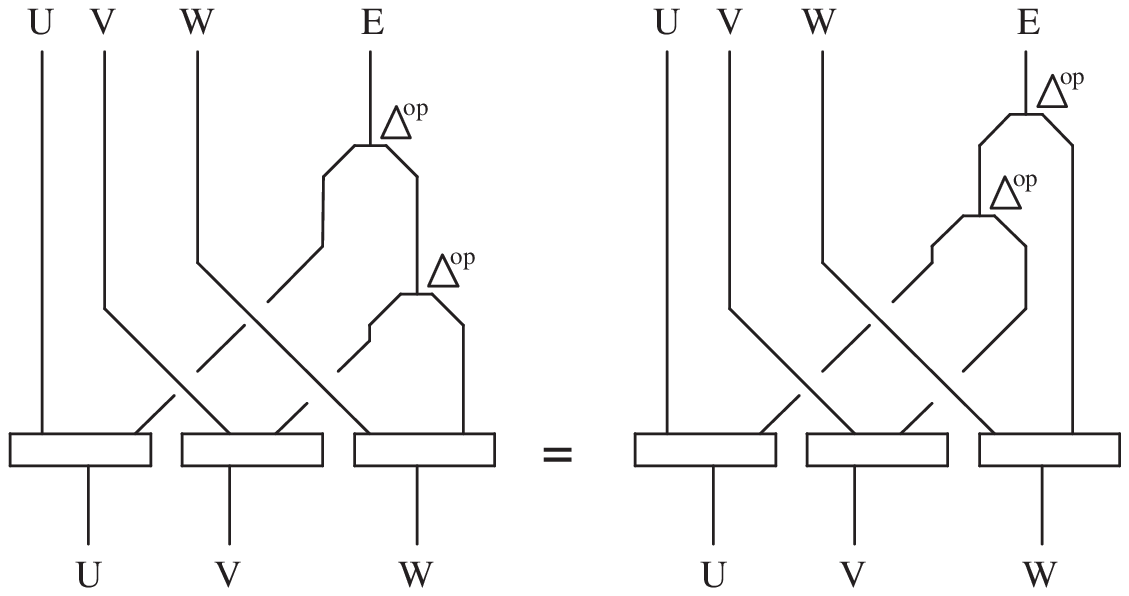}}
\vspace{-1.8cm} \centerline{\rm{ Figure 13} \qquad\qquad}
\end{picture}

\end{Proposition}

\textbf{Proof.}\\
$$
$$
\begin{picture}(10.,10.)
\put(3.,0.){\put(-1.5,3.){\text{L.H.S=}}
\includegraphics{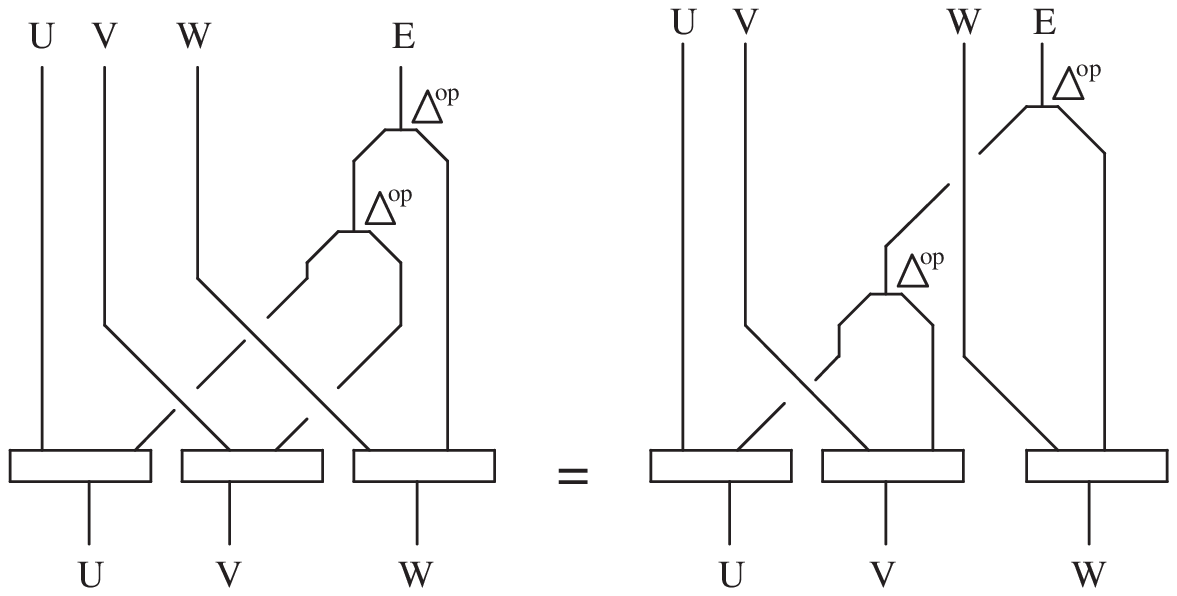}}
\end{picture}
$$
$$
$$
$$
$$
$$
$$
$$
$$
$$
\vspace{.5 cm}
\begin{picture}(10.,10.)
\put(3.,0.){\includegraphics{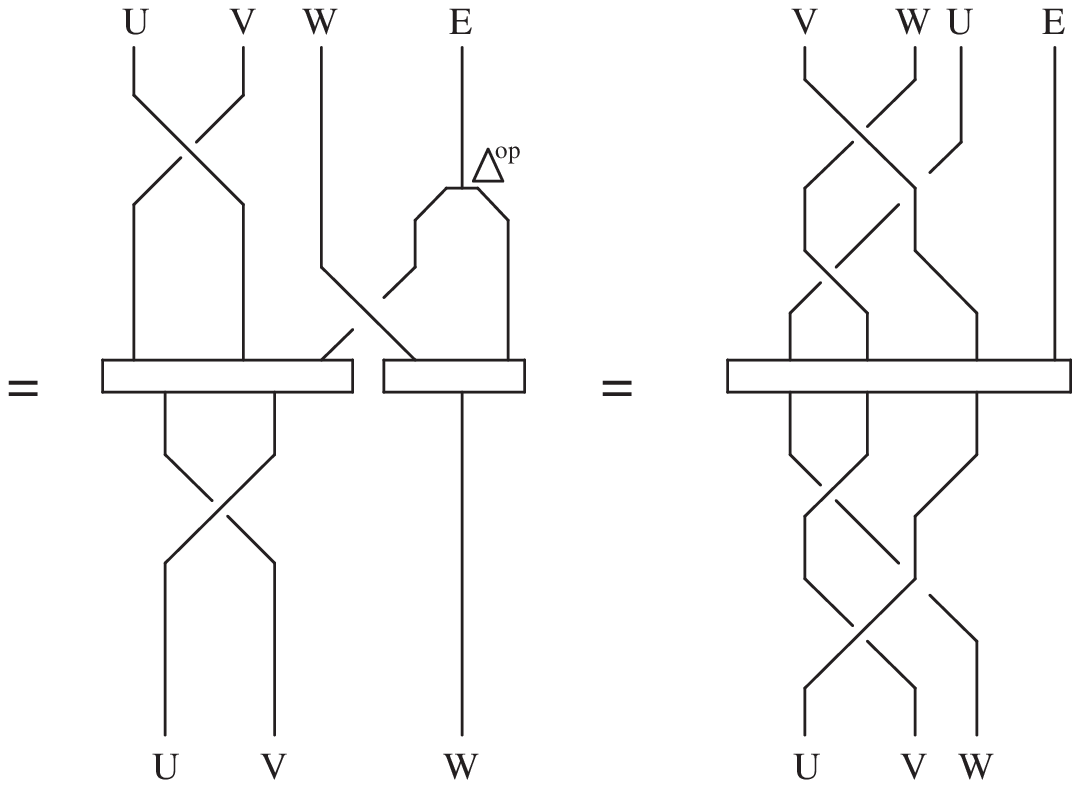}}
\end{picture}\\
$$
$$
$$
$$
$$
$$
$$
$$
\vspace{1. cm}
\begin{picture}(10.,10.)
\put(3.,0.){\put(-1.5,5.5){\text{R.H.S=}}
\includegraphics{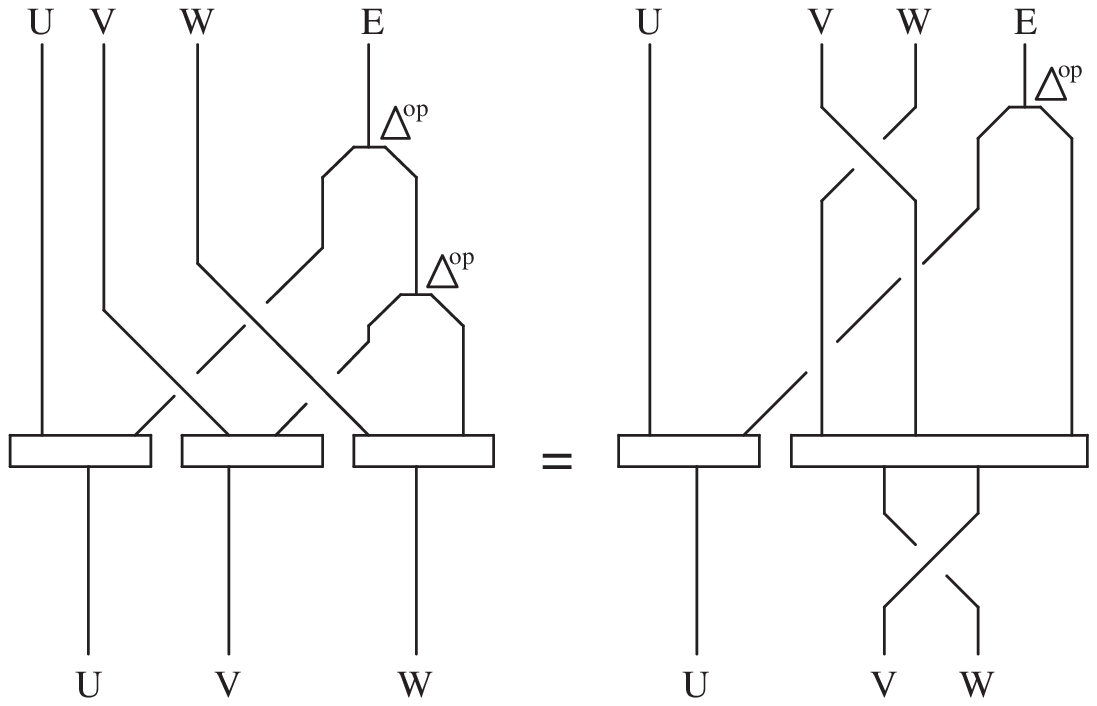}}
\end{picture}
$$
$$
$$
$$
$$
$$
$$
$$
\begin{picture}(10.,10.)
\put(3.,-1.){\includegraphics{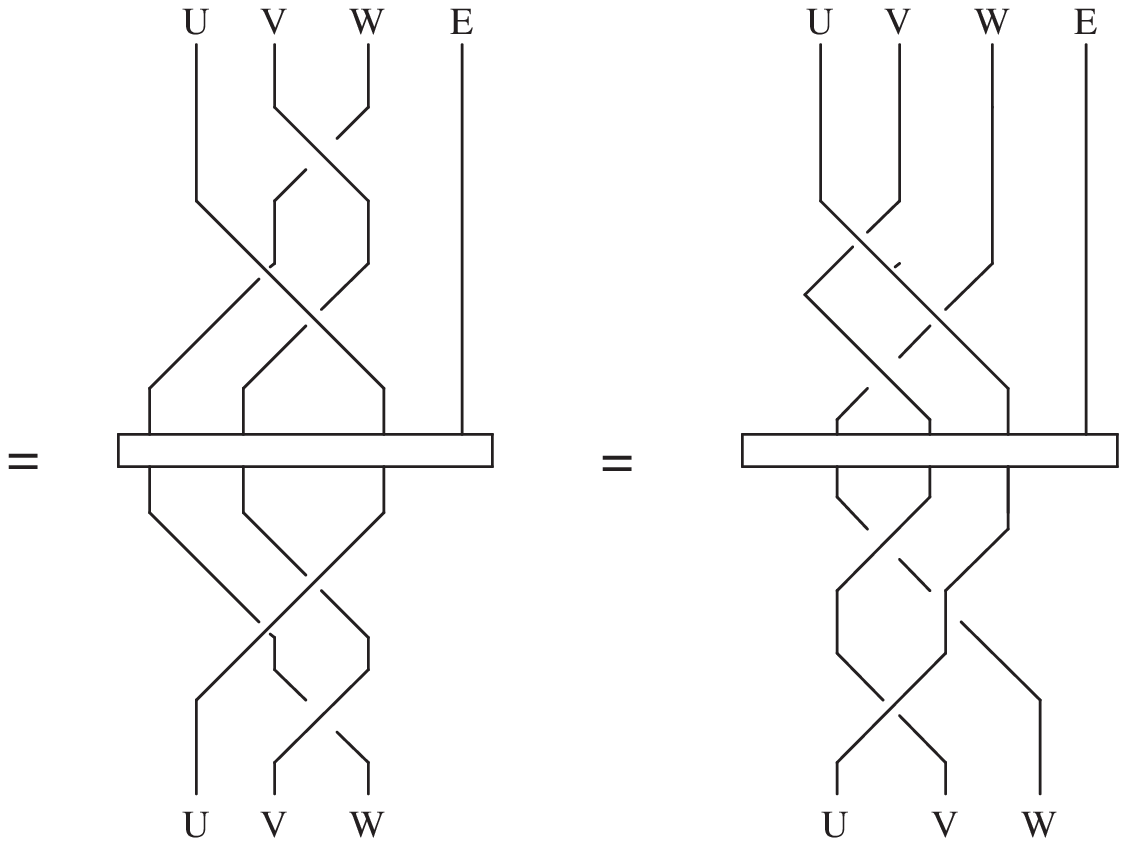}}
\end{picture}\\

\bigskip\noindent
In our case for the category $\mathcal{D}$, we can say more:

\newpage
\begin{Proposition}\,Using the definition of the opposite coproduct
in   \ref{opcopr}, the braided Hopf algebra $D$ in the category
$\mathcal{D}$ is cocommutative.
\end{Proposition}
\textbf{Proof.}\\
$$
$$
%
\begin{picture}(10.,10.)
%

\put(3.5,0.){\includegraphics{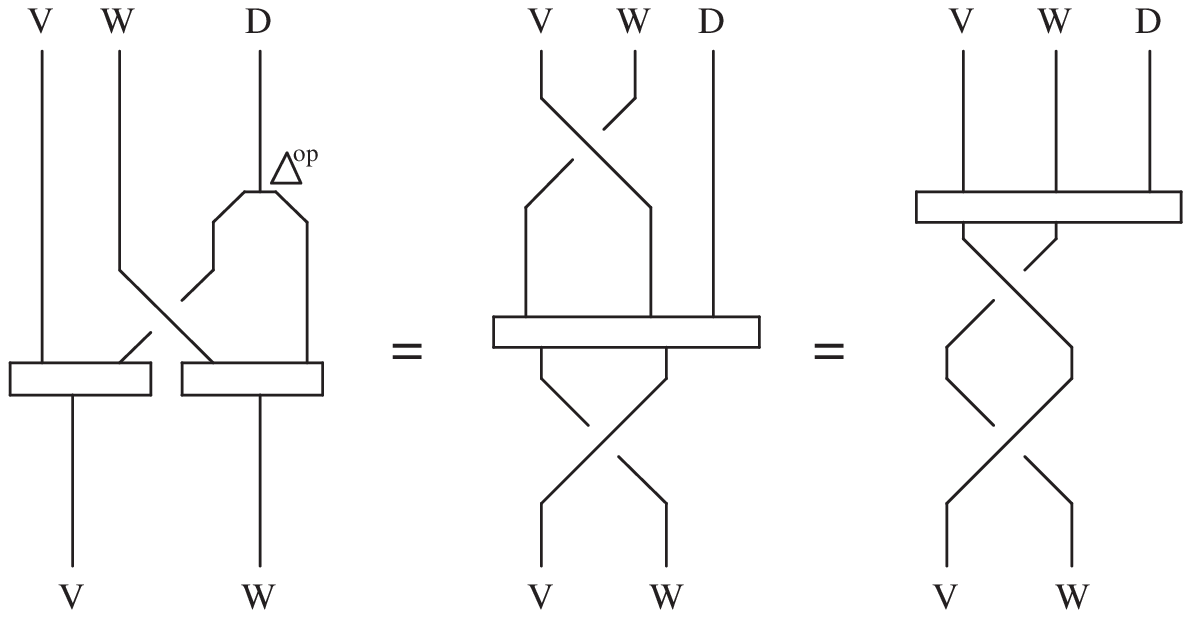}}

\end{picture}\\

%
\qquad\qquad\qquad\quad\begin{picture}(10.,10.)
%

\put(3.5,0.){\includegraphics{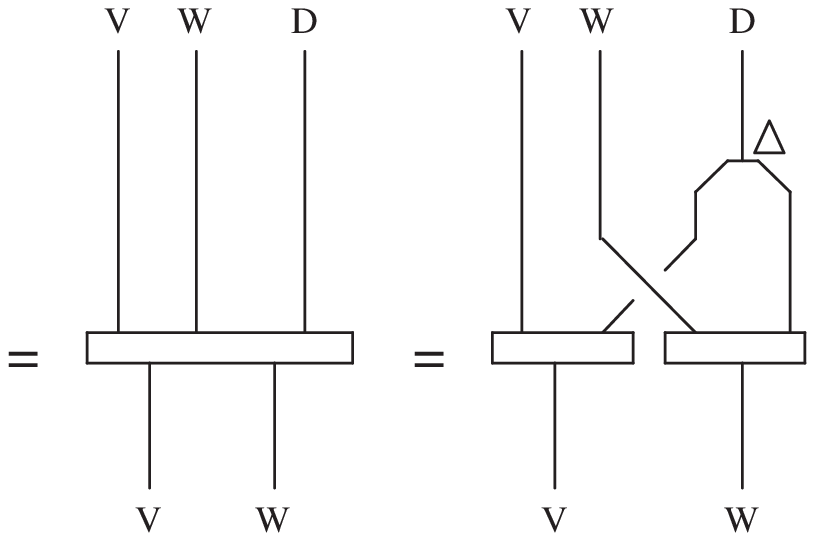}}

\end{picture}
\section{The Hopf algebra $D$ is not braided commutative}

After knowing that the algebra $D$ is braided cocommutative we
would like to know whether it is braided commutative or not, i.e.
whether for $\xi$ and $\eta$ in $D$  the following equation is
satisfied for product $\mu$ and braiding $\Psi$:
\begin{equation}
\begin{split}\label{braidcom}
\mu(\xi \otimes \eta )=\mu \big(\Psi (\xi \otimes \eta)\big)\,\,?
\end{split}
\end{equation}
Put $\xi = \delta_{y} \otimes x$ and $\eta = \delta_{w} \otimes
z$, then the left hand side of (\ref{braidcom}) becomes
\begin{equation}
\begin{split}\label{lbraidcom}
(\delta_{y} \otimes x)  (\delta_{w} \otimes z) = \delta_{w, y
\actl x} \,  \delta_{y \actl \tilde{\tau}(a, b)} \otimes
\tilde{\tau}(a, b)^{-1} x z\,,
\end{split}
\end{equation}
where $a=\| \delta_{y} \otimes x \|=\| \xi \|= \bbi \b$ and $b=\|
\delta_{w} \otimes z \|=\aaa =\aai \a$.  On the other hand \,$
\Psi(\ab) = \eta \achl (\b \acl \aa)^{-1}\otimes \xi \achl \aa
=(\delta_{w} \otimes z)\achl (\b \acl \aa)^{-1}\otimes (\delta_{y}
\otimes x) \achl \aa  $,\, so
\begin{equation}
\begin{split}\label{rbraidcom}
\mu \big(\Psi (\xi \otimes \eta)\big)&= \big(\delta_{w \actl (b
\actr (\b \acl \aa)^{-1})} \otimes (b \actr (\b \acl
\aa)^{-1})^{-1} z (\b \acl \aa)^{-1}\big) \,  \,\big( \delta_{y
\actl (a \actr \aa)} \otimes (a \actr \aa)^{-1} x
\aa \big)\\
&= \delta_{y \actl (a \actr \aa)\,,\,w \actl z (\b \acl \aa)^{-1}}
\, \, \delta_{(w \actl (b \actr (\b \acl \aa)^{-1}))\actl \tau(a
\actl (\b \acl \aa)^{-1} , b \actl \aa)} \\
&\otimes\, \tau(a \actl (\b \acl \aa)^{-1} , b \actl \aa)^{-1} (b
\actr (\b \acl \aa)^{-1})^{-1} z (\b \acl \aa)^{-1} (a \actr
\aa)^{-1} x \aa \,.
\end{split}
\end{equation}
To check the $\delta$ function the statement $w \actl z (\b \acl
\aa)^{-1}(a \actr \aa)^{-1} = y$ should be the same as $w \actl
x^{-1} =y $, i.e.
$$
w \actl z (\b \acl \aa)^{-1}(a \actr \aa)^{-1} x= w\,,
$$
which means
\begin{equation}\label{deltafun}
w  z (\b \acl \aa)^{-1}(a \actr \aa)^{-1} x=z (\b \acl \aa)^{-1}(a
\actr \aa)^{-1} x w\,.
\end{equation}
Now to calculate $(a \actr \aa)(\b \acl \aa)$, put $a=\bbi \b=vt$,
$\b=t$ and $\aa= \bar{w}$ then using the fact that  $vt \actr wp =
v^{-1}wp v^{'}=t wp {t^{'}}^{-1}\,, \text{where}\,\,\, vt \actl wp
=v^{'}t^{'}$ we get
$$
{\bar{w}}^{-1} vt \bar{w}={\bar{w}}^{-1} v (t \acr \bar{w}) (t
\acl \bar{w})=v^{'} t^{'}\,.
$$
So $t^{'}=(t \acl \bar{w})=(\b \acl \aa)$ which implies $(a \actr
\aa)(\b \acl \aa)= t \bar{w} =\b \aa$.  Thus (\ref{deltafun})
becomes
\begin{equation}\label{deltafunc}
w  z (\b \aa)^{-1} x=z (\b \aa)^{-1} x w\,.
\end{equation}
To check whether this is always true or not we consider the
following example from \cite{shombeggs}:
\begin{example}
 Take $X$ to be the dihedral group
$D_6=\langle a,b : a^6=b^2=e, ab=ba^5 \rangle$, whose elements we
list as $\{e, a, a^2, a^3, a^4, a^5, b, ba, ba^2, ba^3, ba^4, ba^5 \}$,
and $G$ to be the non-abelian normal subgroup of order 6 generated by
$a^2$ and $b$, i.e.\ $G=\{e, a^2, a^4, b, ba^2, ba^4\}$.  We choose
$M=\{e, a \}$.  

 Now let $(\delta_{y} \otimes x)
=(\delta_{b{a^{n}}} \otimes b{a^{m}})$ and  $(\delta_{w} \otimes
z) =(\delta_{b{a^{p}}} \otimes b{a^{q}})$. Then we need to check
if the following equation holds:
\begin{equation}\label{deltafunce}
b{a^{p}}  b{a^{q}} (\langle \delta_{b{a^{n}}} \otimes b{a^{m}}
\rangle |\delta_{b{a^{p}}} \otimes b{a^{q}}|)^{-1}
b{a^{m}}=b{a^{q}} (\langle \delta_{b{a^{n}}} \otimes b{a^{m}}
\rangle |\delta_{b{a^{p}}} \otimes b{a^{q}}|)^{-1} b{a^{m}}
b{a^{p}}\,.
\end{equation}
To do so we need to calculate $\| \delta_{b{a^{n}}} \otimes
b{a^{m}} \|$, which we do as follows
$$
b{a^{n}} \circ \| \delta_{b{a^{n}}} \otimes b{a^{m}} \| =b{a^{n}}
\actl b{a^{m}}=(b{a^{m}})^{-1} b{a^{n}} b{a^{m}}=b{a^{2m-n}}\,.
$$
Put $\| \delta_{b{a^{n}}} \otimes b{a^{m}}
\|={b^{\alpha}}{a^{\beta}}$, where $\alpha= 0,1$ and $\beta$ is
even, then ${b^{\alpha}}{a^{\beta}} b{a^{n}}={b^{\alpha +1
}}{a^{n-\beta}} =b{a^{2m-n}}$ which implies $\alpha =0$ and
$\beta=2n - 2m$\,. Thus
$$
\| \delta_{b{a^{n}}} \otimes b{a^{m}} \|=| \delta_{b{a^{n}}}
\otimes b{a^{m}} |^{-1} \langle  \delta_{b{a^{n}}} \otimes
b{a^{m}} \rangle ={a^{2n-2m}} \in G \,,
$$
which implies that $| \delta_{b{a^{n}}} \otimes b{a^{m}}
|={a^{2m-2n}}$ and $\langle  \delta_{b{a^{n}}} \otimes b{a^{m}}
\rangle =e$.  So the left hand side of (\ref{deltafunce}) is
$$
b{a^{p}}  b{a^{q}}  |\delta_{b{a^{p}}} \otimes b{a^{q}}|^{-1}
b{a^{m}}=b{a^{p}}  b{a^{q}} {a^{2p-2q}} b{a^{m}}=b{a^{q-p+m}}\,,
$$
on the other hand  the right hand side of (\ref{deltafunce}) is
$$
b{a^{q}} |\delta_{b{a^{p}}} \otimes b{a^{q}}|^{-1} b{a^{m}}
b{a^{p}}=b{a^{q}} {a^{2p-2q}} b{a^{m}} b{a^{p}} =b {a^{3p-q-m}}\,
$$
which shows that the left hand side of (\ref{deltafunce}) is not
equal to the right hand side, otherwise $q-p+m \equiv_{_6}
3p-q-m$, i.e. $2q-4p+2m$ is a multiple of $6$ which is not always
true. Therefore, we conclude that $D$ is not braided commutative.
\end{example}

\section{Type $A$ and Type $B$ Morphisms}
We will assume for this section that $s^{LL}=s$ and $s^{L} \acr (s
\acr u) = u$ \, for $s \in M$ and $u \in G$ (this is true when $M$ is a
subgroup).

We add new morphisms to the category
$\mathcal{C}$, to make a new category $\tilde{\mathcal{C}}$\,.
Consider the linear map \quad $\phi : V \rightarrow W$. We call it
a type A morphism if it satisfies the following conditions:
$$
\langle \phi(\xi)\rangle = \langle \xi \rangle \qquad
\rm{and}\qquad \phi (\xi \acbl u)=\phi(\xi)\acbl u \,,
$$
for $\xi \in V$ and $u \in G$ (these are just the usual morphism
conditions in $\mathcal{C}$). It  is said to be a type B morphism
if it satisfies the following conditions:
$$
\langle \phi(\xi)\rangle = \langle \xi \rangle^{L} \qquad
\rm{and}\qquad \phi (\xi \acbl u)=\phi(\xi)\acbl (\langle \xi
\rangle \acr u)\,.
$$
We need to check that these morphisms are closed under composition.
It is obvious that the composition of a type A morphism and a type
B morphism is a type B morphism.

\begin{Proposition}\,\,\,
The composition of two type B morphisms is a type A morphism.
\end{Proposition}
\textbf{Proof.}\,\quad\,Let $\phi : U \rightarrow V$ and $\varphi
: V \rightarrow W $ be two   type B morphisms and for $\xi \in U$
let $\b=s$\,.  We first check the grade as the following: As
$\phi$ is a type B morphism then $\langle \phi(\xi) \rangle =
\langle \xi \rangle^{L}$\,. So as $\varphi$ is also a type B
morphism then $\langle \varphi (\phi (\xi)) \rangle
=\b^{LL}=s^{LL}=s$ which is the same as type A morphism.

Now to check the $G$-action, we do the following: As $\phi$ and
$\varphi$ are  type B morphisms then we have the following:
$$
\phi (\xi \acbl u)=\phi(\xi)\acbl (\langle \xi \rangle \acr u)
\qquad\rm{and}\qquad \varphi (\eta \acbl u)=\varphi(\eta)\acbl
(\langle \eta \rangle \acr u)\,,
$$
for $\eta \in V$.  So their composition can be given as the
following:
\begin{equation*}
\begin{split}
\varphi \big(\phi(\xi \acbl u) \big)&=\varphi\big(\phi(\xi)\acbl
(\langle \xi \rangle \acr u)\big)=\varphi\big(\phi(\xi)\big)\acbl
\big(\langle\phi(\xi)\rangle\acr (\langle \xi \rangle \acr
u)\big)\\
&=\varphi\big(\phi(\xi)\big)\acbl \big(\langle\xi\rangle^{L} \acr
(\langle \xi \rangle \acr
u)\big)=\varphi\big(\phi(\xi)\big)\acbl\, \,u\,,
\end{split}
\end{equation*}
which is also the same as type A morphism.\quad $\square$

\medskip
However we no longer have the usual sort of tensor category, as the
type B morphisms obey a rather odd order reversing tensor product
rule, as we now see:
\begin{Proposition}\label{typebcom}\,\,\,
If $\phi:V \rightarrow \tilde{V}$ and  $\psi: W \rightarrow
\tilde{W}$ are type B morphisms, then the map $\phi \boxtimes \psi
: V \otimes W \rightarrow \tilde{W} \otimes \tilde{V}$ which is
defined by
$$
(\phi \boxtimes \psi) (\ab)= \big( \psi (\eta) \acbl \tau(a^{L}, a
)^{-1} \otimes \phi (\xi)\big) \acbl \tau(a , b )\,,
$$
where $\xi \in V$, $\eta \in W$,  $a=\b$ and $b=\a$, is a type B
morphism.
\end{Proposition}
\textbf{Proof.}\,\quad\,First we need to show that $\langle (\phi
\boxtimes \psi) (\ab) \rangle = \langle \ab \rangle^{L}= (a \cdot
b )^{L}$ which we do as the following, taking into account that
$\phi$ and $\psi$ are type B morphisms:
\begin{equation*}
\begin{split}
\big\langle \big( \psi (\eta) \acbl \tau(a^{L}, a )^{-1} \otimes
\phi (\xi)\big) \acbl \tau(a , b ) \big\rangle &=\big\langle \big(
\psi (\eta) \acbl \tau(a^{L}, a )^{-1} \otimes \phi (\xi)\big)
\big\rangle
\acl \tau(a , b ) \\
&=\big( \big\langle  \psi (\eta) \acbl \tau(a^{L}, a
)^{-1}\big\rangle \cdot \big\langle\phi (\xi)  \big\rangle \big)
\acl \tau(a , b ) \\
&=\big( \big(\langle  \psi (\eta)\rangle \acl \tau(a^{L}, a
)^{-1}\big) \cdot \big\langle\phi (\xi)  \big\rangle \big)
\acl \tau(a , b ) \\
&=\big( \big(\langle  \eta \rangle^{L} \acl \tau(a^{L}, a
)^{-1}\big) \cdot \langle \xi  \rangle^{L} \big)
\acl \tau(a , b ) \\
&=\big( \big(b^{L} \acl \tau(a^{L}, a )^{-1}\big) \cdot a^{L}
\big) \acl \tau(a , b ) \,.
\end{split}
\end{equation*}
To show that this is equal to $(a \cdot b )^{L}$ we dot it by $(a
\cdot b )$ to get the identity.  So we get
\begin{equation*}
\begin{split}
\Big(\big( \big(b^{L} \acl \tau(a^{L}, a )^{-1}\big) \cdot a^{L}
\big) \acl \tau(a , b )\Big) \cdot (a \cdot b )&=\Big(\big(
\big(b^{L} \acl \tau(a^{L}, a )^{-1}\big) \cdot a^{L} \big) \cdot
a \Big)  \cdot b
\\
&=\Big(\big( b^{L} \acl \tau(a^{L}, a )^{-1}\tau(a^{L}, a )\big)
\cdot (a^{L}  \cdot a )\Big)  \cdot b \\
&= b^{L} \cdot b =e\,.
\end{split}
\end{equation*}
Next we need to show that $\big(\phi \boxtimes \psi \big)
\big((\ab)\acbl u \big)=\big((\phi \boxtimes \psi) (\ab)\big)\acbl
(\langle \ab\rangle \acr u )$\,.  We start with the left hand side
as the following:
\begin{equation*}
\begin{split}
L. H. S.&=\big(\phi \boxtimes \psi \big) \big(\xi \acbl(\a \acr
u)\otimes \eta \acbl u \big)=\big(\phi \boxtimes \psi \big)
\big(\xi \acbl(b \acr u)\otimes \eta \acbl u \big)\\
&=\Big( \psi (\eta \acbl u) \acbl \tau\big((a \acl (b \acr
u))^{L}, a \acl (b \acr u) \big)^{-1} \otimes \phi (\xi \acbl  (b
\acr
u))\Big) \acbl \tau\big( a \acl (b \acr u) , b\acl u \big)\\
&=\Big( \big(\psi (\eta )\acbl( \a \acr u)\big) \acbl \tau\big((a
\acl (b \acr u))^{L}, a \acl (b \acr u) \big)^{-1} \otimes \phi
(\xi) \acbl (\b \acr(b \acr
u))\Big) \acbl \tau\big( a \acl (b \acr u) , b\acl u \big)\\
&=\Big( \big(\psi (\eta )\big)\acbl( b \acr u) \, \tau\big((a \acl
(b \acr u))^{L}, a \acl (b \acr u) \big)^{-1} \otimes \phi (\xi)
\acbl (a \acr(b \acr u))\Big) \acbl \tau\big( a \acl (b \acr u) ,
b\acl u \big)\,.
\end{split}
\end{equation*}
On the other hand
\begin{equation*}
\begin{split}
R. H. S.&=\big((\phi \boxtimes \psi) (\ab)\big)\acbl (\langle
\ab\rangle \acr u )=\Big(\big( \psi (\eta) \acbl \tau(a^{L}, a
)^{-1} \otimes \phi (\xi)\big) \acbl \tau(a , b
)\Big)\acbl(\langle
\ab\rangle \acr u )\\
&=\big( \psi (\eta) \acbl \tau(a^{L}, a )^{-1} \otimes \phi
(\xi)\big) \acbl \tau(a , b )\,\big((\b
\cdot \a) \acr u \big)\\
&=\big( \psi (\eta) \acbl \tau(a^{L}, a )^{-1} \otimes \phi
(\xi)\big) \acbl \tau(a , b )\,\big((a
\cdot b) \acr u \big)\\
&=\big( \psi (\eta) \acbl \tau(a^{L}, a )^{-1} \otimes \phi
(\xi)\big) \acbl \big(a \acr (b \acr u) \big)\,\tau\big( a \acl (b
\acr u) , b\acl
u \big)\,\,\\
&=\Big( \big(\psi (\eta )\big)\acbl \, \tau(a^{L}, a  )^{-1} (b
\acr u) \otimes \phi (\xi) \acbl (a \acr(b \acr u))\Big) \acbl
\tau\big( a \acl (b \acr u) , b\acl u \big)\,,
\end{split}
\end{equation*}
which is the same as the left hand side as $\tau(a^{L}, a  )^{-1}
(b \acr u)=( b \acr u) \, \tau\big( a^{L} \acl(a\acr (b \acr u)),
a \acl (b \acr u) \big)^{-1}=( b \acr u) \, \tau\big((a \acl (b
\acr u))^{L}, a \acl (b \acr u) \big)^{-1}$\,. Note that we have
used $a^{L} \acr \big( a \acr (b \acr u)\big)=b \acr u$ \,\,by an
assumption for this section.\quad $\square$

\medskip This tensor product of type B morphisms has the following
composition with the braiding $\Psi$:

\begin{Proposition}\,\,\,
Let $\phi:V \rightarrow \tilde{V}$,\,\,  $\psi: W \rightarrow
\tilde{W}$ and $\phi \boxtimes \psi : V \otimes W \rightarrow
\tilde{W} \otimes \tilde{V}$ be as defined in proposition
\ref{typebcom}.  Then the following equality is satisfied
\begin{equation}\label{main}
\big(\Psi\circ(\psi \boxtimes \phi)\big) (\eta \otimes \xi)=
\big((\phi \boxtimes \psi)\circ \Psi^{-1}\big) (\eta \otimes
\xi) \quad\forall \eta \otimes \xi\in V \otimes W\ .
\end{equation}
\end{Proposition}
\textbf{Proof.}\,\quad\,Using the double construction and
remembering that
$$
\Psi(\ab) = \eta \achl (\b \acl \aa)^{-1}\otimes \xi \achl \aa \,,
\quad\Psi^{-1}(\xi^{'} \otimes \eta^{'}) = \eta^{'} \achl
|\xi^{'} \achl \langle \eta^{'} \rangle |^{-1}\otimes \xi^{'}
\achl \langle \eta^{'} \rangle \,,
$$
we start with the left hand side of (\ref{main}) as the following
\begin{equation*}
\begin{split}
(\psi \boxtimes \phi) (\eta \otimes \xi)&= \big( \phi (\xi) \achl
\tau(\a^{L}, \a )^{-1} \otimes \psi (\eta)\big) \achl \tau(\a , \b
)\\
&=  \phi (\xi) \achl \tau(\a^{L}, \a )^{-1} \big( \a^{L}\acr
\tau(\a , \b )\big)\otimes \psi (\eta) \achl \tau(\a , \b )\,.
\end{split}
\end{equation*}
Now applying the braiding  map to the previous equation gives
\begin{equation}\label{lhsmain}
\Psi\big((\psi \boxtimes \phi) (\eta \otimes
\xi)\big)=\eta^{\prime}\achl (\langle \xi^{\prime} \rangle  \acl
|\eta^{\prime}|)^{-1} \otimes \xi^{\prime} \achl
|\eta^{\prime}|\,,
\end{equation}
where \,$\eta^{\prime}=\psi (\eta) \achl \tau(\a , \b )$ and
\,$\xi^{\prime}=\phi (\xi) \achl \tau(\a^{L}, \a )^{-1} \big(
\a^{L}\acr \tau(\a , \b )\big)$. To simplify  equation
(\ref{lhsmain}) we need to calculate the following
$$
|\eta^{\prime}|=|\psi (\eta) \achl \tau(\a , \b )|=\big(
\a^{L}\acr \tau(\a , \b ) \big)^{-1} |\psi (\eta)| \tau(\a , \b ).
$$
We do not know what $|\psi (\eta)|$ is, but we do know that
$$
||\psi (\eta)||=||\eta||^{L}=(\aai \a)^{L}=\aa \ai = \aa
\tau(\a^{L} , \a )^{-1}\a^{L}\,,
$$
and on the other side $||\psi (\eta)||=|\psi (\eta)|^{-1} \langle
\psi (\eta) \rangle \,$ which implies that $|\psi
(\eta)|=\tau(\a^{L} , \a ) \aai$ by the uniqueness of
factorization. So the right part of the tensor of the right hand
side  of (\ref{lhsmain}) becomes
\begin{equation}\label{lhsmain1}
\xi^{\prime} \achl |\eta^{\prime}|=
\phi (\xi) \achl \tau(\a^{L}, \a )^{-1}|\psi (\eta)| \tau(\a , \b )=
\phi (\xi) \achl \aai \tau(\a , \b )\,.
\end{equation}
Next, for the other part of the tensor we need to calculate
\begin{equation*}
\begin{split}
\langle \xi^{\prime} \rangle  \acl |\eta^{\prime}|&=\langle
\xi^{\prime} \achl |\eta^{\prime}| \rangle=\langle \phi (\xi)
\achl \aai \tau(\a , \b ) \rangle\\
&=\langle \phi (\xi)  \rangle \acl \aai \tau(\a , \b )=\b^{L} \acl
\aai \tau(\a , \b )\,.
\end{split}
\end{equation*}
So the left part of the tensor of the right hand side  of
(\ref{lhsmain}) becomes
\begin{equation}\label{lhsmain2}
\begin{split}
\eta^{\prime}\achl (\langle \xi^{\prime} \rangle  \acl
|\eta^{\prime}|)^{-1}= \psi (\eta) \achl \tau(\a , \b )\big(
\b^{L} \acl \aai\, \tau(\a , \b ) \big)^{-1}.
\end{split}
\end{equation}
Thus from (\ref{lhsmain1}) and (\ref{lhsmain2}), equation
(\ref{lhsmain}) can be rewritten as
\begin{equation}\label{nlhsmain}
\Psi\big((\psi \boxtimes \phi) (\eta \otimes \xi)\big)=\psi (\eta)
\achl \tau(\a , \b )\big( \b^{L} \acl \aai\, \tau(\a , \b )
\big)^{-1} \otimes \phi (\xi) \achl \aai \tau(\a , \b ).
\end{equation}
Now  we turn to   the right hand side of (\ref{main})
\begin{equation}\label{rhsmain}
\begin{split}
\big((\phi \boxtimes \psi)\circ \Psi^{-1}\big) &(\eta \otimes
\xi)=\big(\phi \boxtimes \psi \big) \big(\Psi^{-1} (\eta \otimes
\xi)\big)=\big(\phi \boxtimes \psi \big)\big(\xi \achl |\eta\achl
\b|^{-1} \otimes \eta \achl \b\big)\qquad\\ \, &=\Big(\psi(\eta
\achl \b )\achl \tau\big( (\b \acl u)^{L}, \b \acl
u\big)^{-1}\otimes \phi(\xi \achl u)\Big)\achl \tau (\b \acl u,
\langle\eta\achl \b\rangle)
\\ \, &=\psi(\eta \achl \b
)\achl \tau\big( (\b \acl u)^{L}, \b \acl u\big)^{-1}\big(\|
\phi(\xi \achl u)\|\acr \tau (\b \acl u, \langle\eta\achl
\b\rangle) \big)\\
&\quad\otimes \phi(\xi \achl u)\achl \tau (\b \acl u,
\langle\eta\achl
\b\rangle)\\
\, &=\psi\achl (\eta)(\aaa \actr \b ) \tau\big( (\b \acl u)^{L},
\b \acl u\big)^{-1}\big((\b  \acl u)^{L}\acr \tau (\b \acl u,
\langle\eta\achl
\b\rangle) \big)\\
&\quad\otimes \phi(\xi)\achl(\bbb\actr u) \tau (\b \acl u,
\langle\eta\achl \b\rangle)\,,
\end{split}
\end{equation}
where $u=|\eta\achl \b|^{-1}$.  To simplify equation
(\ref{rhsmain}) we need to make the following calculations
$$
\bbb\actr u=\bbi \b \actr u=\bb u v^{\prime},
$$
where $\bbb\actl u=u^{-1}\bbi \b u=u^{-1}\bbi (\b\acr u)(\b\acl
u)=v^{\prime}t^{\prime}$, which implies, by the uniqueness of
factorization, that $v^{\prime}=u^{-1}\bbi (\b\acr u)$ and hence
$\bbb\actr u= (\b\acr u)$.  Also
\begin{equation*}
\begin{split}
(\b \acl u)\langle\eta\achl \b\rangle&=(\b \acr u)^{-1} \b u
\langle\eta\achl \b\rangle=(\b \acr u)^{-1} \b |\eta\achl \b|^{-1}
\langle\eta\achl \b\rangle \\
&=(\b \acr u)^{-1} \b \|\eta\achl \b\|=(\b \acr u)^{-1} \b
(\|\eta\|\actl \b)\\
&=(\b \acr u)^{-1} \b (\b^{-1}\|\eta\| \b)=(\b \acr u)^{-1}
|\eta|^{-1}\a \b\\
&=(\b \acr u)^{-1} |\eta|^{-1}\tau(\a , \b)(\a \cdot  \b),
\end{split}
\end{equation*}
but on the other hand $(\b \acl u)\langle\eta\achl \b\rangle=\tau
(\b \acl u, \langle\eta\achl \b\rangle) \big((\b \acl u)\cdot
\langle\eta\achl \b\rangle \big)$.  So, by the uniqueness of
factorization, we get
\begin{equation}\label{taudot}
\tau (\b \acl u, \langle\eta\achl \b\rangle)=(\b \acr u)^{-1}
|\eta|^{-1}\tau(\a , \b), \,\,\, \text{and}\,\,\,\big((\b \acl
u)\cdot \langle\eta\achl \b\rangle \big)= (\a \cdot  \b).
\end{equation}
Thus the right part of the tensor of (\ref{rhsmain}) becomes
\begin{equation}\label{rofr}
\phi(\xi)\achl(\bbb\actr u) \tau (\b \acl u, \langle\eta\achl
\b\rangle)=\phi(\xi)\achl |\eta|^{-1}\tau(\a , \b),
\end{equation}
which agrees with the right part of the tensor of
(\ref{nlhsmain}). Now, we need to calculate the following to
simplify the left part of the tensor of (\ref{rhsmain})
\begin{equation*}
\begin{split}
\aaa \actl \b&=\bi \aai \a \b=\tau(\b^{L} ,
\b)^{-1}\b^{L}\aai\tau(\a ,
\b)(\a \cdot \b)\\
&=\tau(\b^{L} , \b)^{-1}\big(\b^{L}\acr\aai\tau(\a ,
\b)\big)\big(\b^{L}\acl\aai\tau(\a ,
\b)\big)(\a \cdot \b)\\
&=\tau(\b^{L} , \b)^{-1}\big(\b^{L}\acr\aai\tau(\a ,
\b)\big)\tau\big((\b^{L}\acl\aai\tau(\a , \b)),(\a \cdot \b)
\big)\\
&\qquad\big(\b^{L}\acl\aai\tau(\a , \b)\big)\cdot(\a \cdot
\b)=v^{\prime\prime} t^{\prime\prime}.
\end{split}
\end{equation*}
So,
\begin{equation*}
\begin{split}
\aaa \actr \b&=\aai \a \actr \b=\aa \b
v^{\prime\prime}\\
&=\aa \b\tau(\b^{L} , \b)^{-1}\big(\b^{L}\acr\aai\tau(\a ,
\b)\big)\tau\big(\b^{L}\acl\aai\tau(\a , \b),\a \cdot \b \big).
\end{split}
\end{equation*}
Also
\begin{equation*}
\begin{split}
(\b  \acl u)^{L}\acr \tau (\b \acl u, \langle\eta\achl \b\rangle)
=&\tau \big((\b \acl u)^{L}, \b \acl u \big)\\
&\tau\Big((\b  \acl u)^{L}\acl \tau (\b \acl u, \langle\eta\achl
\b\rangle), (\b  \acl u)\cdot \langle\eta\achl \b\rangle \Big)^{-1}\\
=&\tau \big((\b \acl u)^{L}, \b \acl u \big)\\
&\tau\Big((\b^{L}\acl(\b \acr u)) \acl (\b \acr u)^{-1}
|\eta|^{-1}\tau(\a , \b), \a \cdot  \b \Big)^{-1}\\
=&\tau \big((\b \acl u)^{L}, \b \acl u \big)\tau\Big(\b^{L}\acl
|\eta|^{-1}\tau(\a , \b), \a \cdot  \b \Big)^{-1}
\end{split}
\end{equation*}
So the  left part of the tensor of (\ref{rhsmain}) can be
rewritten as
\begin{equation*}\label{lofr}
\begin{split}
\psi\achl \aa \b\tau(\b^{L} , \b)^{-1}\big(\b^{L}\acr\aai\tau(\a ,
\b)\big)&=\psi\achl \aa (\b^{L})^{-1}\big(\b^{L}\acr\aai\tau(\a ,
\b)\big)\\
&=\psi\achl \tau(\a , \b)\big(\b^{L}\acl\aai\tau(\a ,
\b)\big)^{-1}.
\end{split}
\end{equation*}
Hence (\ref{rhsmain}) can be rewritten as
\begin{equation*}\label{newrhsmain}
\begin{split}
\big((\phi \boxtimes \psi)\circ \Psi^{-1}\big) &(\eta \otimes
\xi)=\psi\achl \tau(\a , \b)\big(\b^{L}\acl\aai\tau(\a ,
\b)\big)^{-1} \phi(\xi)\achl |\eta|^{-1}\tau(\a , \b).
\end{split}
\end{equation*}
which is the same as the left hand side of
(\ref{main}).\quad$\square$

\medskip
Now we should ask what the effect of a type B morphism is on the
action of the algebra $A$.  The answer is  given in the following
proposition.
\begin{Proposition}\,\,\,
If $\phi: V \rightarrow W$ is a type B morphism, then\\
$$
\setlength{\unitlength}{0.5cm}
\begin{picture}(10.,10.)
\put(.,1.){\includegraphics[scale=1.2]{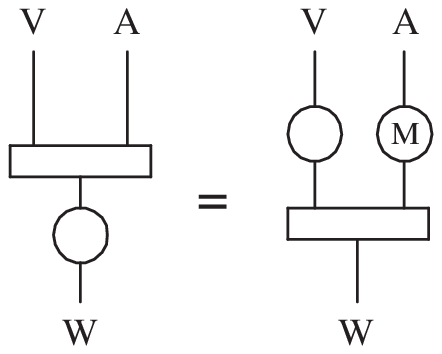}}
\put(3.4,4.2){\text{$\phi$}} \put(9.1,6.7){\text{$\phi$}}
\put(5.8,1){\text{$\rm{Figure\,\, 14}$}}
\end{picture}\qquad\qquad\qquad
$$
where the map $M:A \rightarrow A$ is defined by $M(\delta_{s}
\otimes u )=\delta_{s^{L}} \otimes s \acr u$\,.
\end{Proposition}
\textbf{Proof.}\,\quad\,We start with the left hand side as the
following:  Let $\xi \in V$ and $(\delta_{s} \otimes u)\in A $,
then as \,$\xi \acbl (\delta_{s} \otimes u) = \delta_{s , \b} \,
\xi \acbl u$\, we get
\begin{equation*}
\begin{split}
L. H. S. =\phi \big( \xi \acbl (\delta_{s} \otimes u) \big)=\phi
\big( \delta_{s , \b}  \, \xi \acbl u \big)= \delta_{s , \b} \phi
(   \, \xi \acbl u )\,.
\end{split}
\end{equation*}
To have a non-zero answer we should have $\b =s $.  As $\phi$ is a
type B morphism then
\begin{equation*}
\begin{split}
L. H. S. = \delta_{s , \b} \phi (  \xi ) \acbl ( \b \acr u )=\phi
(  \xi ) \acbl ( s \acr u )\,.
\end{split}
\end{equation*}
Now we calculate the right hand side as the following:
\begin{equation*}
\begin{split}
R. H. S. &=  \phi (  \xi ) \acbl M( \delta_{s} \otimes u )=\phi (
\xi ) \acbl ( \delta_{s^{L}} \otimes s \acr u ) \\
&=\delta_{s^{L} , \langle\phi(\xi)\rangle }\, \phi ( \xi ) \acbl (
s \acr u )=\delta_{s^{L} , \b^{L} }\, \phi ( \xi ) \acbl ( s \acr
u )=\phi ( \xi ) \acbl ( s \acr u ) \,.\quad\square
\end{split}
\end{equation*}

\medskip
 From \cite{BM1}, in the case where $M$ is a subgroup of $X$, there
is a $*$ operation defined on $A$ by $(\delta_{s} \otimes u)^{*}=
\delta_{s \acl u } \otimes u^{-1} $. In our case we have a similar
operation, $P: A \rightarrow A$, given as follows, noting that we
have not yet shown that this really is any sort of conjugation.
\begin{Proposition}\,\,\,\label{themapp}
The map $P: A \rightarrow A$ which is defined  by
$$
P(\delta_{s} \otimes u)= \delta_{s \acl u \tau(a^{L} , a)^{-1}}
\otimes \tau(a^{L} , a)\, u^{-1} \,,
$$
where $a=\langle \delta_{s} \otimes u
\rangle$, is a type B morphism\,.
\end{Proposition}
\textbf{Proof.}\,\quad\,First we check the grade, i.e. $\langle
P(\delta_{s} \otimes u) \rangle = a^{L}$.  It is known that $s
\cdot a = s \acl u$\,.  Now let $\langle P(\delta_{s} \otimes u)
\rangle = b $ and $\tau(a^{L} , a)=w$, then
$$
(s \acl u w^{-1}) \acl w\, u^{-1}=s \acl u w^{-1} w\, u^{-1}=s=(s
\acl u w^{-1})\cdot b \,,
$$
which implies that
$$
s \cdot a =\big((s \acl u w^{-1})\cdot b \big) \cdot a=\big(s \acl
u w^{-1}\tau(b , a) \big)\cdot( b  \cdot a)\,.
$$
But $s \cdot a=s \acl u$, which implies that $b=a^{L}$ as
required.

Now we check the $G$-action, i.e. $P\big((\delta_{s} \otimes u)
\acbl v \big) = P(\delta_{s} \otimes u) \acbl \big( \langle
\delta_{s} \otimes u \rangle \acr v \big)$\,.  We start with the
left hand side as the following:
\begin{equation}
\begin{split}
P\big((\delta_{s} \otimes u) \acbl v \big)&=P\big(\delta_{s \acl
(a \acr v )} \otimes (a \acr v )^{-1} uv \big)\\
&=\delta_{s \acl uv \tau ( (a \acl v)^{L}\,,\,a \acl v)^{-1} }
\otimes \tau ( (a \acl v)^{L}\,,\,a \acl v)\,v^{-1} u^{-1}(a \acr
v )\,.
\end{split}\label{ptypeb}
\end{equation}
To  simplify the last equation we need to do the following
calculation:  Note that $a^{L} \cdot a =e$, so $(a^{L} \cdot a)
\acl v = e \acl v=e$ or $\big(a^{L} \acl (a \acr v)\big) \cdot (a
\acl v )=e$, which means that $\big(a^{L} \acl (a \acr v)\big)=(a
\acl v )^{L}$.  Thus $\tau \big( (a \acl v)^{L}\,,\,a \acl v \big
)^{-1}=\tau \big( (a^{L} \acl (a \acr v))\,,\,a \acl v
\big)^{-1}$, which, from the identities between $(M,\cdot)$ and
$\tau$, implies that
$$
v\,\tau \big( (a \acl v)^{L}\,,\,a \acl v \big )^{-1}=v\,\tau
\big( (a^{L} \acl (a \acr v))\,,\,a \acl v \big)^{-1}= \tau
(a^{L}\,,\,a)^{-1} \big( a^{L} \acr (a \acr v)\big)\,.
$$
So equation (\ref{ptypeb}) can be rewritten as
\begin{equation*}
\begin{split}
P\big((\delta_{s} \otimes u) \acbl v \big)=\delta_{s \acl u\tau
(a^{L}\,,\,a)^{-1} ( a^{L} \acr (a \acr v) ) } \otimes \big( a^{L}
\acr (a \acr v)\big)^{-1}\tau (a^{L}\,,\,a) u^{-1}(a \acr v )\,.
\end{split}
\end{equation*}
On the other hand if we put $a \acr v=\bar{v}$, then the right
hand side is given as the following:
\begin{equation*}
\begin{split}
P(\delta_{s} \otimes u) \acbl \big( \langle \delta_{s} \otimes u
\rangle \acr v \big)&=P(\delta_{s} \otimes u) \acbl \big( a \acr v
\big)=P(\delta_{s} \otimes u) \acbl \bar{v} \\
&=\big(\delta_{s \acl u \tau(a^{L} , a)^{-1}} \otimes \tau(a^{L} ,
a)\, u^{-1}\big) \acbl \bar{v}\\
&=\delta_{(s \acl u \tau(a^{L} , a)^{-1})\acl(a^{L} \acr \bar{v}
)} \otimes (a^{L} \acr \bar{v} )^{-1} \tau(a^{L} ,
a)\, u^{-1} \bar{v}\\
&=\delta_{s \acl u \tau(a^{L} , a)^{-1}\,(a^{L} \acr \bar{v} )}
\otimes (a^{L} \acr \bar{v} )^{-1} \tau(a^{L} ,
a)\, u^{-1} \bar{v}\\
&=\delta_{s \acl u \tau(a^{L} , a)^{-1}\,(a^{L} \acr (a \acr v) )}
\otimes \big(a^{L} \acr (a \acr v) \big)^{-1} \tau(a^{L} , a)\,
u^{-1} (a \acr v)\,,
\end{split}
\end{equation*}
which is the same as the left hand side as required. \quad
$\square$

\begin{Proposition}\,\,\,
For  the algebra $A$ the map $P: A \rightarrow A$ defined in
\ref{themapp} satisfies   \,
$$
P\big( P(\delta_{s} \otimes u) \acbl \tau(a , a^{L}) \big) =
\rm{id_{A}}\,,
$$
where $(\delta_{s} \otimes u) \in A$ and  $a=\langle \delta_{s}
\otimes u \rangle $.
\end{Proposition}
\textbf{Proof.}\,\quad\,First note that $s^{LL}=s $ implies $s^{L}
\acl \tau(s , s^{L})= s^{L}$ and $s^{L}= s^{R}$.  Now if we put
$v=\tau(a , a^{L})$ then
$$
P(\delta_{s} \otimes u) \acbl v= \delta_{s \acl u \tau(a^{L} ,
a)^{-1} (a^{L} \acr v)} \otimes (a^{L} \acr v)^{-1} \tau(a^{L} ,
a)\, u^{-1} v\,.
$$
But $(a^{L} \acr v)^{-1}=\big(a^{L} \acr \tau(a , a^{L})\big)^{-1}
=\tau \big(a^{L} \acl \tau(a , a^{L})\,,\, a \cdot a^{L} \big)
\tau (a^{L} \cdot a , a^{L})^{-1} \tau (a^{L}  , a)^{-1} =\tau
(a^{L} , a)^{-1}$\,,  so
$$
P(\delta_{s} \otimes u) \acbl \tau(a , a^{L})=\delta_{s \acl u}
\otimes  u^{-1} \tau(a , a^{L})\,.
$$
Applying $P$ to this again gives
\begin{equation*}
\begin{split}
P\big( P(\delta_{s} \otimes u) \acbl \tau(a , a^{L}) \big)
&=P\big( \delta_{s \acl u} \otimes  u^{-1} \tau(a , a^{L}) \big)\\
&= \delta_{(s \acl u)\acl u^{-1}\tau(a , a^{L})\tau\big((a^{L}
\acl v)^{L} , a^{L}\acl v\big)^{-1}}
\otimes \tau\big((a^{L} \acl v)^{L} , a^{L}\acl v\big) \tau(a , a^{L})^{-1}
u\\
&= \delta_{s \acl \tau(a , a^{L})\tau(a^{LL} , a^{L} )^{-1}}
\otimes \tau(a^{LL} , a^{L}) \tau(a , a^{L})^{-1} u\\
&= \delta_{s \acl \tau(a , a^{L})\tau(a , a^{L} )^{-1}}
\otimes \tau(a , a^{L}) \tau(a , a^{L})^{-1} u\\
&= \delta_{s } \otimes  u \,.\quad\square
\end{split}
\end{equation*}

\section{\!\!A connection between type $A$ and type $B$ morphisms}

In this section we assume that there is a right inverse in $M$,
and that there is a conjugate $x\rightarrowtail \bar{x}$ on the
field  $k$.
\begin{Def}\label{barfun}
Define a functor $\,Bar: \mathcal{C} \rightarrow \mathcal{C}\,$ as
for $V\in \mathcal{C}$, $Bar(V)=\bar{V}$ where $\bar{V}=V$ as a
set with the usual addition and for $\bar{\xi} \in \bar{V}$,
\,\,$\bar{\xi} x =\xi \bar{x}$ (conjugate scalar multiplication).
In addition, the grade of $\bar{\xi} \in \bar{V}$ is given by
$\langle \bar{\xi} \rangle = \b^{R}$  and the $G$-action on
$\bar{V}$ is given by
$$
\bar{\xi} \, \acbl  \, u = \xi \acbl (\langle \bar{\xi} \rangle
\acr u)= \xi \acbl (\langle \xi \rangle^{R} \acr u).
$$
Moreover, for a morphism $\phi$ in the category,
$\bar{\phi}(\xi)=\phi(\xi)$ as a function between sets.
\end{Def}
\begin{Proposition}\,\,\,
The $M$-grading and the $G$-action given in definition
$\ref{barfun}$ are consistent.
\end{Proposition}
\textbf{Proof.}\,\quad\,
$$
\langle \bar{\xi} \, \acbl  \, u \rangle=\langle \xi \acbl
(\langle \xi \rangle^{R} \acr u) \rangle=\langle \xi\rangle \acl
(\langle \xi \rangle^{R} \acr u)=\langle \xi \rangle^{R} \acl u
=\langle \bar{\xi }\rangle \acl u,
$$
as required where the third equality is due to $(s\cdot t)\acl u =
(s \acl (t \acr u))\cdot (t \acl u)\,\,$ for $s,t\in M$ and $u \in
G$.$\quad \square$

\begin{Proposition}\,\,\,
There is a natural transformation $\Omega$ between the $\,Bar:
\mathcal{C} \rightarrow \mathcal{C}\,$ functor  and the identity
functor $I_{\mathcal{C}}: \mathcal{C}\rightarrow\mathcal{C}$,
defined by
$$
\Omega_{V}(\bar{\bar{\xi}})=\xi \acbl \tau(\langle \bar{\bar{\xi}}
\rangle^{L}, \langle \bar{\bar{\xi}}  \rangle),
$$
that is the following diagram commutes \vspace{2.0cm}
$$\qquad\quad
\setlength{\unitlength}{0.5cm}
\begin{picture}(5,2)\thicklines
\put(-.8,4.8){\hbox to 65pt{\rightarrowfill}} \put(-.80,1){\hbox
to65pt{\rightarrowfill}} \put(-1.8,3){$\Bigg\downarrow$}
\put(4.2,3){$\Bigg\downarrow$}
\put(-1.9,4.8){\text{$\bar{\bar{V}}$}} \put(4.1,4.8){\text{$V$}}
\put(-1.9,1){\text{$\bar{\bar{W}}$}} \put(4.1,1){\text{$W$}}
\put(1.2,5.3){{\text{$\Omega_{V}$}}}
\put(1.2,1.5){{\text{$\Omega_{W} $}}}
\put(-1.4,3.){{\text{$\bar{\bar{\phi}}$}}}
\put(4.6,3.){{\text{$\phi $}}}
\end{picture}
\vspace{-1.0cm}
$$
$$
$$
\end{Proposition}
\textbf{Proof.}\,\quad\,We use
$\xi,\,\bar{\xi},\,\bar{\bar{\xi}}\,$ to distinguish $\,\xi \in
V\,$ as an element of $V,\,\bar{V},\, \bar{\bar{V}}\,$
respectively.  To  show that $\Omega_{V}: \bar{\bar{V}}
\rightarrow V$ is a morphism in the category we need to check the
$M$-grade and the $G$-action.  First we check the $M$-grade as the
following
$$
\langle \Omega_{V}(\bar{\bar{\xi}})  \rangle =\langle \xi \acbl
\tau(\langle \bar{\bar{\xi}}   \rangle^{L}, \langle
\bar{\bar{\xi}}  \rangle) \rangle = \langle \xi\rangle \acl
\tau(\langle {\xi}   \rangle^{R}, \langle {\xi}
\rangle^{RR})=\langle {\xi} \rangle^{RR}=\langle \bar{\bar{{\xi}}}
\,\rangle,
$$
as required.  Now to check the $G$-action, we need to calculate
\begin{equation*}
\begin{split}
\bar{\bar{\xi}}\acbl u &= \bar{{\xi}} \acbl (\langle \bar{{\xi}}
\rangle^{R} \acr u)=\bar{{\xi}} \acbl (\langle {\xi} \rangle^{RR}
\acr u)=\xi \acbl \big(\langle {\xi} \rangle^{R} \acr(\langle
{\xi} \rangle^{RR} \acr u) \big)\\
&=\xi \acbl \tau(\langle {\xi} \rangle^{R} ,\langle {\xi}
\rangle^{RR})\,
u\,\tau\big(\langle {\xi} \rangle^{R} \acl(\langle
{\xi} \rangle^{RR} \acr u) \,,\,(\langle {\xi} \rangle^{RR} \acl
u)\big)^{-1},
\end{split}
\end{equation*}
where the last equality is because $s \acr (t \acr u)=
\tau(s,t)((s \cdot t) \acr u)\tau(s\acl(t\acr u),t \acl u)^{-1}\,$
for $s,t \in M$ and $u \in G$. If we put $\bar{\bar{\xi}}\acbl u
=\bar{\bar{\eta}}$, then
$$
\Omega_{V}(\bar{\bar{\xi}}\acbl
u)=\Omega_{V}(\bar{\bar{\eta}})=\eta \acbl
\tau(\langle\bar{\bar{\eta}}\rangle^{L},
\langle\bar{\bar{\eta}}\rangle)=\eta \acbl \tau(\a^{R}, \a^{RR}).
$$
So
\begin{equation*}
\begin{split}
\a &= \b \acl \tau(\langle {\xi} \rangle^{R} ,\langle {\xi}
\rangle^{RR})\,
u\,\tau\big(\langle {\xi} \rangle^{R} \acl(\langle
{\xi} \rangle^{RR} \acr u) \,,\,(\langle {\xi} \rangle^{RR} \acl
u)\big)^{-1}\\
&=\b^{RR} \acl
u\,\tau\big(\langle {\xi} \rangle^{R} \acl(\langle
{\xi} \rangle^{RR} \acr u) \,,\,(\langle {\xi} \rangle^{RR} \acl
u)\big)^{-1}.
\end{split}
\end{equation*}
Then there is $v \in G$\, such that
\begin{equation*}
\begin{split}
v \a &= \b^{RR} u\,\tau\big(\langle {\xi} \rangle^{R} \acl(\langle
{\xi} \rangle^{RR} \acr u) \,,\,(\langle {\xi} \rangle^{RR} \acl
u)\big)^{-1}\\
&=(\b^{RR}\acr u)(\b^{RR}\acl u)\,\tau\big(\langle {\xi}
\rangle^{R} \acl(\langle {\xi} \rangle^{RR} \acr u) \,,\,(\langle
{\xi} \rangle^{RR} \acl u)\big)^{-1}\\
&=(\b^{RR}\acr u) \big( \b^{R}\acl(\b^{RR}\acr u)\big)^{-1},
\end{split}
\end{equation*}
which implies that $\a = \big( \b^{R}\acl(\b^{RR}\acr
u)\big)^{L}$.  Thus
\begin{equation*}
\begin{split}
\Omega_{V}(\bar{\bar{\eta}})&=\Omega_{V}(\bar{\bar{\xi}}\acbl u)
=\eta \acbl \tau\big( \b^{R}\acl(\b^{RR}\acr u), \b^{RR}\acl
u\big)\\
&=\xi \acbl \tau(\langle {\xi} \rangle^{R} ,\langle {\xi}
\rangle^{RR})u =(\xi \acbl \tau(\langle \bar{\bar{\xi}}
\rangle^{L}, \langle \bar{\bar{\xi}}  \rangle)) \acbl u\\
&=\Omega_{V}(\bar{\bar{\xi}}) \acbl u,
\end{split}
\end{equation*}
as required.$\quad \square$

\begin{Rem}\,\,\,
A type $B$ morphism $\, \phi: V \rightarrow W\,$ can be viewed as
a type $A$ morphism $\, \phi: V \rightarrow \bar{W}\,$ (same as a
function on sets).  Indeed, as for the $M$-grade we have
$$
\langle  \bar{\phi(\xi)} \rangle=\langle  \phi(\xi)
\rangle^{R}=\b.
$$
And for the $G$-action we know $\phi(\xi \acbl u)= \phi(\xi) \acbl
(\b \acr u)$, but we also have
$$
\bar{\phi(\xi)} \acbl u = \phi(\xi)\acbl (\langle  \phi(\xi)
\rangle^{R} \acr u)=\phi(\xi) \acbl (\b \acr u),
$$
as required where $\xi \in V$ $\quad \square$

\end{Rem}

\section{Inner product}

\begin{Def}
An inner product on an object $V$ of the category is given by a
type $B$ morphism $\,\phi: V \rightarrow V^{*}\,$ and then
$$
\langle \eta \,,\,\xi \rangle=eval \big( \phi(\eta)\,,\,\xi
\big)\,,
$$
where  $\eta,\xi \in V$, i.e.
$$
\setlength{\unitlength}{0.5cm}
\begin{picture}(10.,10.)
\put(.,1.){\includegraphics[scale=1.2]{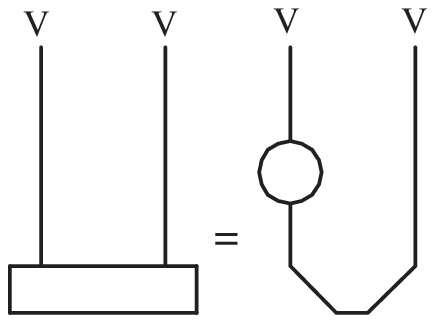}}
\put(3.4,2.32){\text{$\langle \,\,,\,\rangle$}}
\put(8.5,5.2){\text{$\phi$}} \put(5.8,0){\text{$\rm{Figure\,\,
15}$}}
\end{picture}\qquad\qquad\qquad
$$
\end{Def}

\begin{Proposition}\,\,\,
The inner product with a type morphism $\,\phi: V \rightarrow
V^{*}\,$  as defined above is invariant in the category
$\,\mathcal{D}\,$.
\end{Proposition}
\textbf{Proof.}\,\quad\,
$$
\setlength{\unitlength}{0.5cm}
\put(-5.5,0){\begin{picture}(10.,10.)
\put(-5.5,-0.5){\includegraphics[scale=1.2]{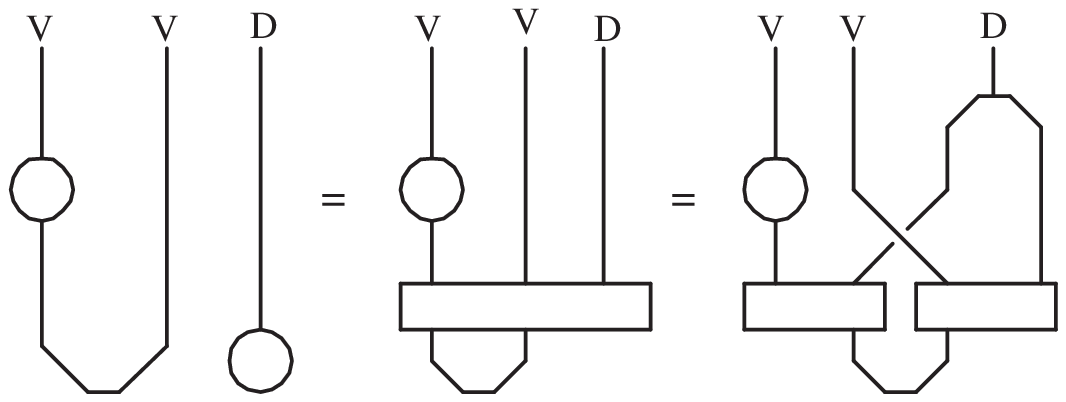}}
\put(-3.9,5.2){\text{$\phi$}} \put(5.6,5.2){\text{$\phi$}}
\put(14,5.2){\text{$\phi$}} \put(1.5,1.){\text{$\epsilon$}}
\put(9,2.3){\text{$\acbl$}}
\end{picture}}\qquad\qquad\qquad
$$
$$
$$
$$
$$
$$
\setlength{\unitlength}{0.5cm} \put(-2.5,0){\begin{picture}(10.,10.)
\put(-5.5,-0.5){\includegraphics[scale=1.2]{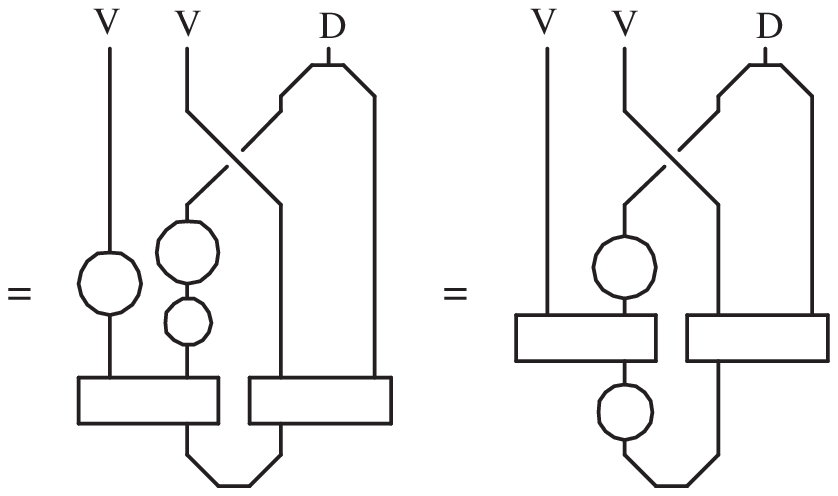}}
\put(-2.8,5.2){\text{$\phi$}}
\put(9.5,5.6){\scriptsize{\text{$M^{-1}$}}}
\put(-1.15,6){\scriptsize{\text{$M^{-1}$}}}
\put(-.9,4.3){\scriptsize{\text{$M$}}}
\put(9.85,2.05){\text{$\phi$}}
\end{picture}}\qquad\qquad\qquad
$$
As an example for the inner product we give the inner product on
the algebra $A$

\begin{example}\,\,
Let $\phi: A \rightarrow A^{*}$\,\, be a type $B$ morphism defined
by
\begin{equation}\label{defoffi}
\phi(\delta_{s } \otimes  u)= s \otimes  \delta_{u }.
\end{equation} Then the inner product on $A$ is given by
$$
\langle \delta_{t } \otimes  v ,  \delta_{s } \otimes  u
\rangle=eval\big( \phi(\delta_{t } \otimes  v) , \delta_{s }
\otimes  u \big),
$$
for $(\delta_{t } \otimes  v)$ and $(\delta_{s } \otimes  u)$ in
$A$.
\end{example}
$\mathbf{Proof.}$\,\, Given the basis $(\delta_{s } \otimes  u)$
of $A$, the dual basis is $\, (s \otimes  \delta_{u })\,$ (see
\cite{BGM}) and also
$$
eval\big( \phi(\delta_{t } \otimes  v) , \delta_{s } \otimes  u
\big)=eval( t \otimes \delta_{v } \, ,\, \delta_{s } \otimes u
)=\delta_{s,t }\,\,\delta_{u,v }.
$$
To have a non-zero solution we must have $t=s$ and $v=u$ and then
we must have $\,\langle  s \otimes  \delta_{u } \rangle=\langle
\delta_{s } \otimes  u \rangle^{L}\,$ as required.

Next, knowing that the evaluation map is invariance under the
action of $G$  we have
$$
eval\big( (t \otimes \delta_{v })\acl (a \acr w) \, ,\, (\delta_{s
} \otimes u)\acl w \big)=\delta_{s,t }\,\,\delta_{u,v }\,,
$$
where $\,a=\langle  \delta_{s } \otimes u \rangle\,$ and $w$ is an
element of the group $G$.  Applying the action we get
$$
eval\big( (t \otimes \delta_{v })\acl (a \acr w) \, ,\, (\delta_{s
\acl (a \acr w) } \otimes  (a \acr w)^{-1} u w \big)=\delta_{s,t
}\,\,\delta_{u,v }\,,
$$
which implies that
\begin{equation}\label{theacofdual}
(t \otimes \delta_{v })\acl (a \acr w)=(s \otimes \delta_{u })\acl
(a \acr w)=s \acl (a \acr w) \otimes  \delta_{(a \acr w)^{-1} u w}
\end{equation}
Now we want to prove that
$$
\phi\big((\delta_{s } \otimes  u)\acbl w \big)= \phi(\delta_{s }
\otimes u)\acbl (a \acr w).
$$
Starting with the left hand side we apply the action then
(\ref{defoffi}) to get
$$
\phi\big((\delta_{s } \otimes  u)\acbl w \big)=\phi\big(\delta_{s
\acl (a \acr w) } \otimes  (a \acr w)^{-1} u w\big)=s \acl (a \acr
w) \otimes \delta_{(a \acr w)^{-1} u w }.
$$
For the right hand side we apply (\ref{defoffi}) then
(\ref{theacofdual}) to get
$$
\phi(\delta_{s } \otimes u)\acbl (a \acr w)=(s \otimes \delta_{u
}) \acbl (a \acr w)=s \acl (a \acr w) \otimes \delta_{(a \acr
w)^{-1} u w },
$$
as required. $\quad \square$

\begin{Def}\,\,
Let $\,\phi\,$ be a type $B$ morphism and let $\,\rho\,$ as
defined in theorem \ref{roint}.   Then the star operation $*$ can
be defined by:
$$
$$
$$
\setlength{\unitlength}{0.5cm}
\put(-4.0,0.0){\begin{picture}(10.,10.)
\put(-1.0,-0.5){\includegraphics[scale=1.2]{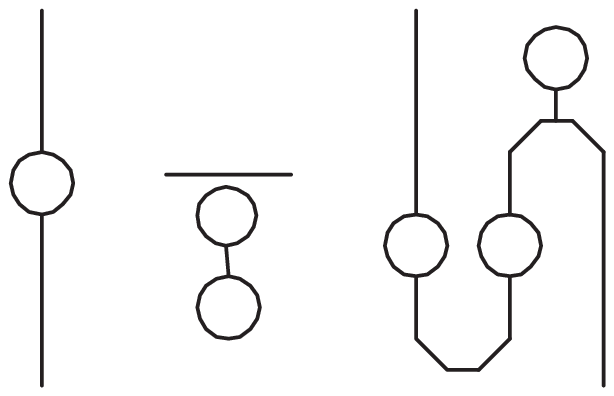}}
\put(3.3,4.8){\text{\rm{\Bigg(\qquad\qquad\Bigg)}}}
\put(5.3,5.7){\text{\rm{1}}} \put(2.3,5){\text{\rm{=}}}
\put(.7,4.5){\text{\rm{*}}} \put(15.5,5.0){\text{\rm{.}}}
\put(5.2,4.0){\text{$\rho$}}\put(5.2,1.7){\text{$\int$}}
\put(9.8,3.3){\text{$\phi$}} \put(12.,3.2){\text{$S$}}
\put(13.2,7.9){\text{$\rho$}} \put(9.8,9.5){\text{$V$}}
\put(0.7,9.5){\text{$V$}} \put(5.,-.5){\text{$\rm{Figure\,\,
16}$}}
\end{picture}}\qquad\qquad\qquad
$$
\end{Def}

\begin{theorem}\,\,\,For an object $V$ and a type $B$ morphism $\phi$,
the following equality holds:
$$
\setlength{\unitlength}{0.5cm}
\put(-4.0,0.0){\begin{picture}(10.,10.)
\put(1.0,-0.5){\includegraphics[scale=1.2]{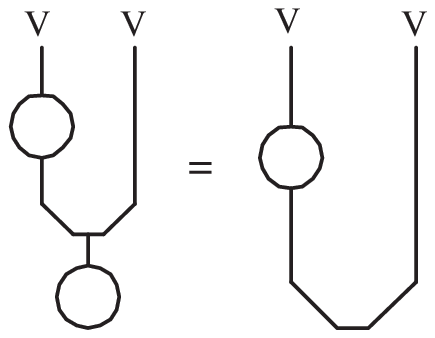}}
\put(2.7,4.9){\text{\rm{*}}} \put(3.7,1.){\text{$\int$}}
\put(8.7,4.5){\text{$\phi$}} \put(5.,-.5){\text{$\rm{Figure\,\,
17}$}}
\end{picture}}\qquad\qquad\qquad
$$
\end{theorem}
\textbf{Proof.}
$$
$$
$$
\setlength{\unitlength}{0.5cm}
\put(-4.0,0.0){\begin{picture}(10.,10.)
\put(-10.0,-0.5){\includegraphics[scale=1.2]{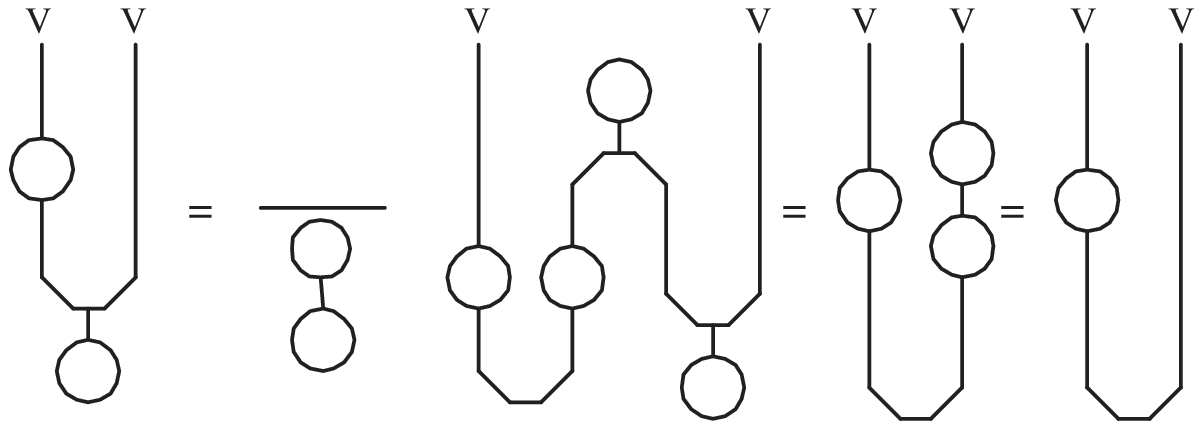}} 
\put(-8.3,6.0){\text{\rm{*}}}
\put(-7.2,1.3){\text{$\int$}} \put(2.3,3.7){\text{$\phi$}}
\put(4.5,3.6){\text{$S$}} \put(-1.5,2.1){\text{$\int$}}
\put(-1.5,4.5){\text{$\rho$}}
\put(-3.5,5){\text{\rm{\Bigg(\qquad\qquad\Bigg)}}}
\put(-1.5,5.9){\text{\rm{1}}} \put(5.8,8.3){\text{$\rho$}}
\put(8.,0.95){\text{$\int$}} \put(11.9,5.6){\text{$\phi$}}
\put(17.2,5.6){\text{$\phi$}}
\put(14.1,4.4){\text{\scriptsize{$S$}}}
\put(13.9,6.6){\text{\scriptsize{$S^{-1}$}}}
\end{picture}}\qquad\qquad\qquad
$$

\end{document}